\documentstyle[11pt,newlfont,amscd,amssymb]{amsart}
\newcommand{\nc}{\newcommand}

\nc{\one}{\mbox{\bf 1}}
\nc{\invtensor}{\underset{\leftarrow}{\otimes}}
\nc{\rlarrows}{\begin{picture}(1,0.4)
                \put(0,-0.1){\makebox(1,0.2){$\leftarrow$}}
                \put(0,0.1){\makebox(1,0.2){$\ra$}}
              \end{picture}}
\nc{\rra}{\begin{picture}(1,0.4)
                \put(0,-0.1){\makebox(1,0.2){$\lra$}}
                \put(0,0.1){\makebox(1,0.2){$\lra$}}
              \end{picture}}
\nc{\Left}{\mathbf L}  
\nc{\Right}{\mathbf R} 
\newcommand{\bs}{\bigskip}

\newcommand{\bm}{{\widehat{{\mathfrak b}}_{-}}}

\newcommand{\nm}{{\widehat{{\mathfrak n}}_{-}}}

\renewcommand{\b}{\beta}

\newcommand{\C}[2]{\left[ \begin{array}{c}#1\\#2\end{array}\right] }

\renewcommand{\dfrac}[2]{\displaystyle{\frac{#1}{#2}}}
\newcommand{\tab}{\hspace*{\fill}}
\newcommand{\ok}{\tab $\Box$\bs}

\newcommand{\SP}{{\Sigma}_{+}}

\newcommand{\SM}{{\Sigma}_{-}}

\newcommand{\SPM}{{\Sigma}_{\pm}}

\newcommand{\SPN}{{\Sigma}_{+,n}}
\newcommand{\SMN}{{\Sigma}_{-,n}}
\newcommand{\SPMN}{{\Sigma}_{\pm,n}}
\newcommand{\SPB}{{\overline{\Sigma}}_+}

\newcommand{\SMB}{{\overline{\Sigma}}_-}

\newcommand{\SEPB}{{\overline{\Sigma}}_{\varepsilon}}
\newcommand{\SEPAB}{{\overline{\Sigma}}_{\varepsilon_0}}
\newcommand{\SEPBB}{{\overline{\Sigma}}_{\varepsilon_1}}

\newcommand{\SEPPB}{{\overline{\Sigma}}_{\varepsilon'}}
\newcommand{\SPMB}{{\overline{\Sigma}}_{\pm}}

\newcommand{\SMPB}{{\overline{\Sigma}}_{\mp}}

\newtheorem{lem}{Lemma}
\newtheorem{prop}{Proposition}
\newtheorem{thm}{Theorem}
\newtheorem{cor}{Corollary}

\newenvironment{Def}{{\noindent\bf Definition. }}{}

\newenvironment{dem}{{\noindent\bf Proof. }}{\ok}

\newenvironment{ex}{{\noindent\bf Examples. }}{}

\renewcommand{\a}{\alpha}
\renewcommand{\b}{\beta}
\newcommand{\NN}{{\mathbb N}}

\newcommand{\CC}{{\mathbb C}}

\newcommand{\QQ}{{\mathbb Q}}
\newcommand{\ZZ}{{\mathbb Z}}
\newcommand{\ootimes}{{\bar\otimes}}
\newcommand{\DDelta}{{\bar\Delta}}
\nc{\fa}{\frak a}
\nc{\fb}{\frak b}
\nc{\fg}{\frak g}
\nc{\fk}{\frak k}
\nc{\fh}{\frak h}
\nc{\fm}{\frak m}
\nc{\fn}{\frak n}
\nc{\fS}{\frak S}
\nc{\fI}{\frak I}
\nc{\fA}{\frak A}
\nc{\nen}{\newenvironment}
\nc{\ol}{\overline}
\nc{\ul}{\underline}
\nc{\ra}{\rightarrow}
\nc{\lra}{\longrightarrow}
\nc{\lla}{\longleftarrow}
\nc{\Lra}{\Longrightarrow}
\nc{\Lla}{\Longleftarrow}
\nc{\Llra}{\Longleftrightarrow}
\nc{\hra}{\hookrightarrow}
\nc{\iso}{\overset{\sim}{\lra}}
\nc{\on}{\operatorname}
\nc{\ad}{\on{ad}}
\nc{\DS}{\on{DS}}
\nc{\DSB}{{\overline{\DS}}}
\nc{\Der}{\on{Der}}
\nc{\Inf}{\on{Inf}}
\nc{\Max}{\on{Max}}
\nc{\card}{\on{card}}
\nc{\mult}{\on{mult}}
\nc{\diag}{\on{diag}}

\newcommand{\cl}{\on{cl}}
\renewcommand{\Im}{\on{Im}}
\renewcommand{\ker}{\on{Ker}}
\begin{document}
\title[]{Discrete quantum Drinfeld-Sokolov correspondence}
\author{C. Grunspan}
\address{C. Grunspan\\
Faculty of Mathematics and Computer Science\\
The Weizmann Institute of Science\\ 
Rehovot 76100, Israel.}

\begin{abstract}
We construct a discrete quantum version of the Drinfeld-Sokolov correspondence
for the sine-Gordon system. The classical version of this correspondence
is a birational Poisson morphism between the phase space of the
discrete sine-Gordon system and a Poisson homogeneous space. Under this
correspondence, the commuting higher mKdV vector fields correspond
to the action of an Abelian Lie algebra. We quantize this picture
(1) by quantizing this Poisson homogeneous space, together with
the action of the Abelian Lie algebra, (2) by quantizing the
sine-Gordon phase space, (3) by computing the quantum analogues of the integrals
of motion generating the mKdV vector fields, and (4) by constructing an algebra
morphism taking one commuting family of derivations to the other one.
\end{abstract}
\maketitle

\section{Introduction}
\subsection{Background}
The link between integrable systems and quantum groups has been 
intensively studied during the last few years from several viewpoints. 
The goal of this article is to present a discrete and quantum version
of a natural construction occurring in the theory of 
integrable systems, namely
the Drinfeld-Sokolov correspondence \cite{DS}. At the classical level, 
this correspondence
was discovered by Drinfeld and Sokolov in the eighties, 
by using the 
dressing method of
Zakharov and Shabat \cite{SZ}. It is a bijective map between phase 
spaces of certain 
evolution equations (such as the KdV, mKdV, or Toda  hierarchy) and 
homogeneous spaces.
Each phase space is equipped with an infinite commuting   family of 
vector fields~:
the Hamiltonian fields generated by the integrals of motion of the
KdV, mKdV, or Toda hierarchy. 
One of the main properties of the Drinfeld-Sokolov correspondence 
is that it leads to a geometric interpretation of these commutative
families~: they
correspond to the action 
of an Abelian Lie algebra on a double coset space. After Drinfeld and 
Sokolov, Feigin and Frenkel
developed a cohomological approach based on the fact that screening 
operators of the Toda theory
satisfy the Serre relations \cite{FEIF} and \cite{FEI}. This allowed 
them (1) to prove 
the existence of a commutative family of integrals of motion in the 
quantum case and (2) to suggest 
a possible discretization of the problem, generalizing those
introduced much earlier by Izergin and Korepin 
(see \cite{IZK}, \cite{IZKO}). 
At the semi-classical level, 
the discretized
Toda system has been studied by Enriquez and Feigin in the case when 
the Lie algebra
is $\widehat{{\mathfrak{sl}}_2}$ \cite{ENR}. This is the discrete 
sine-Gordon theory.
By imitating the cohomological approach of Feigin and Frenkel, 
the authors (1) proved the existence
of a classical family of integrals of motion in 
involution and (2) constructed 
a Drinfeld-Sokolov correspondence between phase spaces 
equipped with the Hamiltonian 
action of the integrals of motion, and homogeneous spaces 
equipped with the action of an
Abelian Lie algebra. 
Moreover, the phase spaces are endowed with a natural structure of
Poisson manifold, 
the homogeneous spaces are Poisson homogeneous spaces and the correspondence
is a Poisson isomorphism. 
The aim of our present work is to quantify this result.
So, this article fills in the discrete-quantum square of the following array.

$$
{\text{\bf{Drinfeld-Sokolov correspondence}}}
$$
\bigskip
$$
\begin{array}{r|c|c}
&\text{classical}&\text{quantum}\\
\hline \text{continuous }&\text{Drinfeld-Sokolov\, (1981)}&\text{?}\\
\hline \text{discrete }&\text{Enriquez-Feigin\, (1995)}&\text{this article\, (2000)}
\end{array}
$$
\subsection{The classical Drinfeld-Sokolov correspondence}

The Drinfeld-Sokolov correspondence is inspired by the application 
of dressing methods
developed by Zakharov and Shabat in the theory of integrable systems
\cite{SZ}. The integrable systems studied by Drinfeld and Sokolov
by this method are
Korteweg-de Vries hierarchies (KdV) or modified
Korteweg-de Vries hierarchies (mKdV), associated with an affine 
Kac-Moody algebra. For example,
in the case when the Kac-Moody algebra ${\tilde{\fg}}$ is
${\widehat{\mathfrak{sl}_2}}$, the second equation of the mKdV 
hierarchy is 
(the first one being $\partial_{z}u=u_z$)
\begin{equation}
u_t=u_{zzz}+6u^2u_z.
\end{equation}
The main achievement of Drinfeld and Sokolov was (1) to associate 
Lax pairs $(A(u),L(u))$ to these equations taking values in 
affine Kac-Moody Lie algebras
and then (2) by assigning to a point of the phase
space, the matrix conjugating its Lax matrix to a prescribed form,
to set up a bijection between the phase space and a coset space
(3) to show that the corresponding system on the coset
space is "linear". This way, Drinfeld and Sokolov 
achieved the linearization   of their system.                            
More precisely, if u belongs to the phase space, the matrices $K(u)$
conjugating
the matrix $L(u)$ into a standard form belong 
to a pro-algebraic
pro-unipotent subgroup $N_+$ of the Kac-Moody group $G$ associated to 
${\tilde{\fg}}$. Moreover, such a matrix $K(u)$ is determined uniquely up to 
a multiplication
by an element of a commutative subgroup $A_+$ in $N_+$, and all the 
coefficients
of the class of $K(u)$ in $N_+/A_+$ are differential polynomials in $u$.
As a result, one gets a map from the phase space of the hierarchy
to $N_+/A_+$. 
The Drinfeld-Sokolov theorem asserts that this map is bijective. 
Moreover, in the corresponding bijection between the rings of functions, 
the hierarchy equations viewed as commutative flows 
on the phase space 
correspond to the right action of the Lie algebra of the normalizer 
of $A_+$ in $G$, $N_+$ being embedded into the flag variety 
$B_-\backslash G$. The Hamiltonian structure (which one 
can associate to these 
hierarchies) were studied
by Gelfand, Dickey and Dorfman (\cite{GDI1},
\cite{GDI2}, \cite{GDO}).

\subsection{The viewpoint of Feigin and Frenkel}

Feigin and Frenkel reformulated the Drinfeld-Sokolov correspondence
in  a cohomological language. 
This allowed them to identify 
the action of $\fn_+$ on the phase space $U$ 
(which is, according to the Drinfeld-Sokolov correspondence, the same as 
the left action by vector fields of $\fn_+$ on the homogeneous space 
$N_+/A_+$) with the Hamiltonian action of screening charges of the Toda system
associated with the Lie algebra ${\tilde{\fg}}$.
Besides, their formalism led to a quantization as well as a 
discretization of the
Toda system.
Precisely, let $\fg$ be a semi-simple Lie algebra and ${\tilde{\fg}}$
be the affine Kac-Moody algebra built from $\fg$.
The Toda system associated with the Lie algebra $\tilde{\fg}$
is the following system of differential 
equations~:
\begin{equation}
\partial_{\tau}\partial_{z}\phi_{i}(z,\tau)=
\sum_{j=0}^{l}(\a_i,\a_j)\exp(-\phi_j(z,\tau)),\,
i=1,\ldots,l,
\end{equation}
where $\a_0,\ldots,\a_l$ are simple roots of ${\tilde{\fg}}$.
Each function $\phi_i(z,\tau)$ depends on $z$ as well as the time 
variable $\tau$,
and $\phi_0(z,\tau)=-\dfrac{1}{a_0}\sum_{i=1}^{l}a_i\phi_i(z,\tau)$,
where $a_0,\ldots,a_l$ are labels of the Dynkin diagram.
In the case when $\tilde{\fg}=\widehat{{\mathfrak{sl}}_2}$,
the system reduces to the sine-Gordon equation~:
\begin{equation}
\partial_{\tau}\partial_{z}\phi(z,\tau)=
\exp(\phi(z,\tau))-\exp(-\phi(z,\tau)).
\end{equation}
Let $\pi_0$ be the ring of functions on $U$. We have 
$\pi_0=\CC[u_i^{(n)}]_{1\leq i\leq l;0<n}$.
This is a differential ring equipped with the derivation
$\partial$ defined by $\partial u_i^{(n)}=u_i^{(n+1)}$.
Because of the presence of $\exp(-\phi_i)$ (and contrary to the mKdV case)
one cannot view the Toda equations
as derivations of the differential ring $\pi_0$.
Yet, one can view the right hand side of these equations
as the action of an evolution operator
(i.e., a linear map commuting with the action of derivation
$\partial$),
\begin{equation}
{\tilde{{\cal H}}}=\sum_{i\in J} \tilde{Q_i}~:\quad \pi_0
\longrightarrow\bigoplus_{i\in J} \pi_{-\alpha_i}
\end{equation} 
for some differential modules $\pi_{-\alpha_i}=\pi_0\otimes\exp(-\phi_i)$
equipped with derivations $\partial$ defined by
$\partial(u\otimes\exp(-\phi_i))=
(\partial u)\otimes\exp(-\phi_i)
-(uu_i^{(0)})\otimes\exp(-\phi_i)$.
One of the results of Feigin and
Frenkel shows that ${\tilde{Q}}_i=-T_i e_i^{G}$, 
where $e_i^{G}$ denotes the image under the Drinfeld-Sokolov correspondence
of the left action by the vector field $e_i$ on $N_+/A_+$, and $T_i$ 
denotes the
multiplication by $\exp(-\phi_i)$ sending $\pi_0$ to
$\pi_{-\alpha_i}$.
In the Hamiltonian formalism of Feigin and Frenkel, one 
defines ${\cal F}_{-\alpha_i}:=\pi_{-\alpha_i}/\text{Im}\partial$.
This is the space of functionals of the form
\begin{equation}
u\longmapsto\int_{|z|=1}P\bigl(u(z),\partial_z u(z),\ldots\bigr)
\exp(-\phi_i(z))dz
\end{equation}
with $u(z)=\bigl(u_1^{(0)},\ldots,u_l^{(0)}\bigr):\, S^1\longrightarrow\fh$,
where $S^1$ is the unit circle and
$\phi_i$ is an anti-derivative of $u_i^{(0)}$.
The space ${\cal F}_0$ is the local functionals space of the Toda system
(one uses only derivatives of $u$). 
It was shown in [GDi,GDo] that it may be equipped with the 
structure of a "vertex Poisson" algebra (this notion was first
introduced in these papers, and developed in [BD,EF]).
Such vertex Poisson structures 
are classical limits of families of VOA structure
(as opposed to classical limits of associative algebra
structures in the case of Poisson structures).
It is possible to extend the Poisson bracket to bilinear maps~:
\begin{equation}
\begin{array}{rcl}
{\cal F}_0\times{\cal F}_{-\alpha_i}&\longrightarrow&{\cal F}_{-\alpha_i}\\
(F,G)&\longmapsto&\{ F,G\}
\end{array}
\end{equation}
satisfying the Jacobi identity on ${\cal F}_0$. In other words, 
the ${\cal F}_{-\alpha_i}$ are ${\cal F}_0$-modules.
By passing to the quotients, the morphisms $\tilde{Q_i}$
define screening operators $\bar{Q_i}:{\cal F}_0
\longrightarrow{\cal F}_{-\alpha_i}$. Feigin and Frenkel showed that 
${\bar Q}_i=\{ .,\int\exp(-\phi_i)\}$, where 
$\int\exp(-\phi_i)$ is the projection of
$\exp(-\phi_i)$ onto ${\cal F}_{-\alpha_i}$.
This projection is called the classical screening charge. Then, we define 
a Hamiltonian 
\begin{equation}
{\cal H}=\sum {\bar Q}_i~:\quad 
{\cal F}_0\longrightarrow \bigoplus_{i\in J} {\cal F}_{-\alpha_i}.
\end{equation}
Indeed, in this formalism, we can write Toda equations as
\begin{equation}
\partial_{\tau} u(z)=\{ u(z),{\cal H}\}.
\end{equation}
The integrals of motion of the Toda system are the local 
functionals which commute with all the screening charges.
On the other hand, it was known \cite{BMP} that the 
the ${\bar Q}_i$ satisfy
the Serre relations for $\fg$. This gives an action of the 
nilpotent Lie algebra 
$\fn_+$ on $\pi_0$ and allowed
Feigin and Frenkel to interpret the space of integrals of motion as the 
first cohomology group of $\fn_+$ with coefficients in $\pi_0$.
Using resolutions of BGG type, the authors managed to compute this
cohomology, yet without giving explicit formulas for the integrals of motions
(except for some particular cases). They showed that the space of integrals of
motion forms a graded one-dimensional Lie subalgebra of ${\cal F}_0$
(for the gradation defined by $\partial^{o}u_i^{(k)}=-k$ on $\pi_0$),
generated by integrals of motion $I_m,\, m\in{\ZZ}_-$ and $\partial^{o} I_m=m$.
This Lie subalgebra is the classical ${\cal W}$-algebra of $\fg$.
Moreover, the Hamiltonian flow of $I_m$ corresponds to the
mKdV hierarchy equations. 
Similarly, this Hamiltonian approach allowed them to show that the integrals of motion
admit quantum deformations.
\subsection{The quantum sine-Gordon system}
The interpretation of the Poisson structure on 
${\cal F}_0$ in terms of Kirillov-Kostant structures allows the 
following quantum deformation (\cite{FEI}).
\subsubsection{The continuous quantum model}
Consider the quantum Heisenberg algebra of 
$\widehat{{\mathfrak{sl}}_2}$, generated by 
$I,\, q_i,\, b_n^{i},\, i\in\{ 0,1\},\, n\in\ZZ$,
satisfying the following relations~:
\begin{align*}
[I,q_i]&=[I,b_n^{i}]=0\\
[q_i,b_n^{j}]&=(\a_i,\a_j)\delta_{n,0}\,\b^{2}I\\
[b_n^{i},b_n^{j}]&=n(\a_i,\a_j)\delta_{n+m,0}\,\beta^{2}I,
\end{align*}
where $\delta_{i,j}$ denotes the Kronecker symbol, 
$\beta^{2}$ is a deformation parameter ($q=\exp(i\pi\beta^{2})$), 
and $(\a_i,\a_j)$ denotes the scalar product of two roots 
$\a_i$ and $\a_j$ in $\widehat{{\mathfrak{sl}}_2}$.

This algebra acts on the direct sum of Fock modules in such a way that the vertex
operators~:
   $$
    V_i(z)=:\exp\bigl(\phi_i(z)\bigr):=\exp\bigl(\phi_i^{-}(z)\bigr)
    \exp\bigl(\phi_i^{+}(z)\bigr),
   $$
satisfy the following commutation relations~:
\begin{equation}\label{vjvi}
  V_j(w)V_i(z)=\exp\bigl(i\pi\beta^{2}(\a_i,\a_j)\bigr)
  V_i(z)V_j(w)
 \end{equation}
in the domain $|z|>|w|$. Here,
$\phi_i(z)$ is the free field~:
   $$
    \phi_i(z)=\sum_{n\not= 0}-\dfrac{1}{n}b_n^{i}z^{-n}
-b_0^{i}\ln(z)-q_i,
   $$
and 
\begin{align*}
\phi_i^{+}(z)&=\sum_{n>0}-\dfrac{1}{n}
b_n^{i}z^{-n}-b_0^{i}\ln(z)\\
\phi_i^{-}(z)&=\sum_{n<0}-\dfrac{1}{n}
b_n^{i}z^{-n}-b_0^{i}\ln(z).
\end{align*}

\noindent
The screening charges are defined by~:
\begin{equation}\label{si}
S_i=\int_{|z|=1}V_i(z)dz=\int_{|z|=1}:\exp\bigl(\phi_i(z)\bigr):dz
\end{equation}
for $i\in\{ 0,1\}$.
The integrals of motion of the system are expressions of the form
   $$
    \int_{|z|=1}P\bigl(\partial_z^{k}\phi_i\bigr)_{k\geq 0,i\in\{
      0,1\}}(z)dz
   $$
which commute with $S_0$ and $S_1$.
Feigin and Frenkel showed that the screening charges $S_0$ and $S_1$ 
satisfy quantum Serre relations of $\widehat{{\mathfrak{sl}}_2}$. 
This allowed them -- as in the classical case -- to interpret 
the space of quantum
integrals of motion as the first cohomology group of a certain complex.
They proved that the integrals of motion commute with each other and that
the space
of integrals of motion generates a quantum deformation
of the classical ${\cal{W}}$-algebra.
\subsubsection{The discrete quantum sine-Gordon model}
The first model of a discrete integrable system was introduced by Izergin
and Korepin in 1982 for the sine-Gordon system in order to resolve 
ultraviolet divergence problems
occurring in the continuous theory. 
The $q$-commutation relations between vertex operators naturally lead to the 
aforementioned discretization, which is adopted 
by Izergin and Korepin.

Set $q=\exp(i\pi\beta^2)$, and replace the complex numbers $z$ by 
relative integers $k\in\ZZ$. The vertex operators $V_0(z)$ and $V_1(z)$
are replaced by variables $y_k$ and $x_k$ satisfying the relations~:
\begin{align}
\label{rotra}
\forall k<l,\quad x_k x_l&=qx_l x_k\\
\label{rotrb}
y_k y_l&=q y_l y_k\\
\label{rotrc}
y_k x_l&=q^{-1}x_l y_k\\
\label{rotrd}
\forall k\leq l,\quad x_k y_l&=q^{-1}y_l x_k,
\end{align}
coming from (\ref{vjvi}). The analogues of screening charges are $S_0$ and $S_1$
defined by $S_0=\sum_{k=-\infty}^{+\infty}y_k$
and $S_1=\sum_{k=-\infty}^{+\infty}x_k$ 
in a certain completion of the algebra of variables on the given lattice.
The Hamiltonian of the system is $H=S_0+S_1$. The integrals of motion 
correspond to expressions of the form
   $$
    \sum_{i=-\infty}^{+\infty}P(x_{i+1}^{\pm 1},y_{i+1}^{\pm
      1},\ldots,x_{i+k}^{\pm 1},y_{i+k}^{\pm 1})
   $$
which commute with $S_0$ and $S_1$, where
$P(X_{1}^{\pm 1},Y_1^{\pm 1},\ldots,X_k^{\pm 1},Y_k^{\pm 1})$
denotes a polynomial in the variables,
$X_1^{\pm 1},\ldots,Y_k^{\pm 1}$, these variables being $q$-commutative.

As in the continuous case, Enriquez and Feigin showed 
that the screening charges satisfy the
Serre relations (for $\widehat{{\mathfrak{sl}}_2}$),
which gives a cohomological interpretation for the integrals of motion.
By means of Demazure desingularization, they managed to compute this
cohomology in the classical limit $q\rightarrow 1$, to give formulae for
densities of integrals of motion, and to prove that
the integrals of motion are in involution. This justifies calling the
system a discrete integrable system.
Moreover, Enriquez and Feigin identified the phase space of this system with
the homogeneous space $H_-\backslash B_-$, where $B_-$ 
is a Borel subgroup of the loop group of $SL_2$, and $H_-$ is a
subgroup of $B_-$ consisting of diagonal matrices, and established
a discrete version of the Drinfeld-Sokolov correspondence. The Hamiltonian
action by integrals of motions corresponds to a (left) action of a commutative
Lie algebra $\fh_+$ on the homogeneous space, which is embedded into 
$H_-\backslash G/N_+$. Moreover, the correspondence is a Poisson morphism.

\section{Main results}\label{resprin}
In this section, we present our 
main results which deal with the discrete sine-Gordon theory.
Proposition \ref{basi} is proved in \cite{Gru}.
In Proposition \ref{thb}, we construct a
quantization of the Poisson homogeneous space considered by 
Enriquez and Feigin
\cite{ENR}. Theorem \ref{thc} is a discrete and quantum version of the 
Drinfeld-Sokolov correspondence. It generalizes a 
theorem of Enriquez and Feigin to the quantum case.

\subsection{Some basic definitions}\label{defprin}
The following notation will be used throughout the article.
Let $q$ be a formal variable.

{\bf The quantum phase space of the discrete sine-Gordon system.}
Let $A_q$ be the algebra generated over $\QQ[q,q^{-1}]$ by the 
non-commutative variables
$x_i^{\pm 1}$ and $y_i^{\pm 1}$, $i\in\ZZ$, subject to relations 
(\ref{rotra}),\ldots, (\ref{rotrd}). This is the quantum phase space 
of our system. 
It can be shown that a basis for $A_q$ is given by the family 
$\prod_{i=-\infty}^{+\infty}
x_i^{\a_i}y_i^{\b_i}$, where $(\a_i)$ and $(\b_i)$ are two sequences in 
$\ZZ^{\ZZ}$ which have
only a finite number of non-zero elements. At the semi-classical limit,
$A_q$ defines a Poisson structure on 
$A_{\cl}=\QQ[x_i^{\pm 1},y_i^{\pm 1},i\in\ZZ]$.
There is a gradation $\deg$ on $A_q$ defined by 
$\deg(x_i)=-\deg(y_i)=1$ for all 
$i\in\ZZ$. Let $T^{\frac{1}{2}}$ denote the half-translation 
antigraded automorphism
defined on $A_q$ by $T^{\frac{1}{2}}(x_i)=y_i$ and 
$T^{\frac{1}{2}}(y_i)=x_{i+1}$ for all $i\in\ZZ$. We shall set
$T={\bigl(T^{\frac{1}{2}}\bigr)}^2$. For $n\in\ZZ$, let $A_q[n]$ 
be the submodule of $A_q$
of all elements of degree $n$. The discrete analogue of 
$\pi_0$ is $A_q[0]$.

{\bf The functional spaces and integrals of motion.}
By definition, the functional spaces are ${\cal F}_n$, $n\in\NN$
defined by ${\cal F}_n=A_q[n]/{\Im(T-Id)}$.
If $P\in A_q[n]$, one notes $I(P)$ its class in ${\cal F}_n$.
The discrete analogues of the screening
charges $S_0$ and $S_1$ seen in (\ref{si}), are $\SP=I(x_0)\in{\cal F}_1$ and 
$\SM=I(y_0)\in{\cal F}_{-1}$. The Hamiltonian of the system is $H=\SM+\SP\in
{\cal F}_{-1}\oplus{\cal F}_1$. We set ${\cal I}=\ker[.,\SM]\cap
\ker[.,\SP]\subset{\cal F}_0$. It is the space of all local integrals of 
motion of the system.
If $P\in A_q[0]$ and if $I(P)\in{\cal I}$, we say that $P$ is a density
of an integral of motion.

{\bf The homogeneous space of Enriquez and Feigin.}
We set $G=\text{SL}_2\bigl(\CC((\lambda^{-1}))\bigr)$,
$B_-=\pi^{-1}\bigl(\bar{B}_-\bigr)$, and $N_+=p^{-1}\bigl(\bar{N}_+\bigr)$,
where  $\bar{B}_-$ (resp. $\bar{N}_+$) is the subgroup of $\text{SL}_2(\CC)$ 
defined by all lower triangular (resp. upper unipotent) matrices, 
$\pi$
(resp. $p$) is the induced map from 
$\text{SL}_2\bigl(\CC[[\lambda^{-1}]]\bigr)$ 
(resp. $\text{SL}_2\bigl(\CC[\lambda]\bigr)$)
to $\text{SL}_2(\CC)$ obtained
by sending $\lambda^{-1}$ (resp. $\lambda$) to $0$. We also denote by
$H_-$ the subgroup of $B_-$ given by all the diagonal matrices of the form
$\diag(a,a^{-1})$, $a\in\CC[[\lambda^{-1}]]^*$. The Poisson homogeneous space
considered by Enriquez and Feigin is 
$H_-\backslash B_-$ endowed with the Poisson structure induced by 
the Poisson bivector
$P_{\infty}=r^L-r_{\infty}^R$ on $G$ such that the map
$H_-\backslash B_-\hookrightarrow H_-\backslash G/N_+$ is a Poisson morphism.
Here, $r$ is the standard $r$-matrix on $G$, $r_{\infty}$ corresponds to the 
conjugation of $r$ with an element of the Weyl group with an infinite length
(it is the $r$-matrix associated with the ``new realizations'' of Drinfeld)
and $r^L$ (resp. $r_{\infty}^R$) is the left (resp. right) translation of $r$
(resp. $r_{\infty}$) on $G$. The homogeneous space $H_-\backslash B_-$
is endowed with an action from the left by the Abelian Lie algebra
$\fh_+:=\{ \diag(a,-a),\, a\in\lambda\CC[\lambda]\}$.
Enriquez and Feigin identified this action with the Hamiltonian action of the 
integrals of motion on $A_{\cl}$. Theorem \ref{thc} 
generalizes this result.
\subsection{The results}
The first few results deal with the integrals of motion of the discrete sine-Gordon system.
\begin{prop}\label{nondema}
Let $n$ be an integer. For any 
$F\in{\cal F}_0,\, G\in{\cal F}_n,\, P\in I^{-1}(F)\subset A_q[0]$ 
and $Q\in I^{-1}(G)\subset A_q[n]$, 
the two sums $\sum_{i=-\infty}^{+\infty}[T^i P,Q]$
and $\sum_{i=-\infty}^{+\infty}[P,T^i Q]$
are equal and do not depend on $P$ nor $Q$ but only of $F$
and $G$.
Moreover, the bilinear map~:
\begin{equation}
\begin{array}{rcl}
{\cal F}_0\times{\cal F}_n&\longrightarrow&{\cal F}_n\\
(F,G)&\longmapsto&[F,G]
\end{array}
\end{equation}
with $F=I(P)$, $G=I(Q)$, and 
\begin{align}
\label{alij1}
[I,J]&:=\frac{1}{q-1}\sum_{i=-\infty}^{+\infty}[T^i P,Q]\\
\label{alij2}
&=\frac{1}{q-1}\sum_{i=-\infty}^{+\infty}[P,T^i Q]
\end{align}
satisfies the Jacobi identity on ${\cal F}_0$.
In other words,
the space
${\cal F}_0$ of all local functionals is a Lie algebra, and ${\cal F}_n$
is a  ${\cal F}_0$-module.
\end{prop}
The following proposition shows that the local functionals act on $A_q$
by adjoint action.

\begin{prop}\label{nondemb}
For any $F\in{\cal F}_0$ and $x\in A_q$, the sum
$\sum_{k=-\infty}^{\infty}[T^k P,x]$ does not depend
of the representative $P\in I^{-1}(F)\subset A_q[0]$ chosen for $F$
and defines a derivation on $A_q$ which commutes with the automorphism
$T$. This derivation is the adjoint action of $F$ on $A_q$
and is denoted by $\ad(F)$.
Moreover, if $\Der_{T}(A_q)$ denotes the Lie algebra of all derivations 
on $A_q$ which commute with $T$, we have a  Lie algebra morphism~:
\begin{equation}\label{adjoin}
\begin{array}{rcl}
\ad~:\quad{\cal F}_0&\longrightarrow&\Der_{T}(A_q)\\
F&\longmapsto&\frac{1}{q-1}\sum_{k=-\infty}^{\infty}[T^k P,.]
\end{array}
\end{equation}
and the kernel of this morphism
is equal to the class of constants in ${\cal F}_0$.
\end{prop}

Proposition \ref{basi} gives an explicit basis for the space ${\cal I}$
of local integrals of motion. This should be compared to the
formulas involving quantum trace identities of Izergin and Korepin \cite{IZKO}
(see also \cite{FTT} and \cite{HIK}) at least in the classical case.
Indeed, in the quantum case, the Izergin and Korepin formulas for
the sine-Gordon system are no longer local.

\begin{prop}\label{basi}
The local classical integrals of motion of the discrete sine-Gordon system can be 
deformed in the quantum case. 
A basis for ${\cal I}$ is given by the family $I_n=I(\psi_n),\, n>0$.
The generating function of the densitites of integrals of motion $\psi_n$ 
is given by 
\begin{equation}
\ln_q U +\ln_q V=\sum_{p=1}^{\infty}\dfrac{1}{[p]}\psi_p\lambda^{-p}
\end{equation}
where for all integers $p$,
$[p]:=\dfrac{q^p-1}{q-1}$,
$U=\lim\limits_{N\to\infty}U_N$, $V=\lim\limits_{N\to\infty}V_N$,

\begin{align}
U_N&:=\dfrac{1}{1-\dfrac{(\lambda x_1y_1)^{-1}}
      {1-\dfrac{(\lambda y_{1}x_2)^{-1}}{\dfrac{\ddots
      (\lambda x_{N-1}y_{N-1})^{-1}}
      {1-(\lambda y_{N-1}x_{N})^{-1}}}}}\\
V_N&:=T^{{1\over 2}}U_N
\end{align}
the limits being taken in the sense of the $\lambda^{-1}-topology$.
By convention, we have set $\dfrac{a}{b}:=ab^{-1}$, and for all power series 
$f$ in $\lambda^{-1}$ with non-zero constant term, 
$$
\ln_q f:= \sum_{p=1}^{\infty} \dfrac{1}{[p]}(1-f^{-1})^p.
$$
\end{prop}
The following proposition is proved in section \ref{unom}.
\begin{prop}\label{tha}
The space ${\cal I}$ 
is a commutative Lie subalgebra in ${\cal F}_0$.
\end{prop}
Proposition \ref{thb} gives a quantization of the homogeneous space 
considered by Enriquez and Feigin.
\begin{prop}\label{thb}
The algebra  given by generators~: $u_i,m_i,\, i>0$ and 
(quadratic) relations
\begin{align}
\label{relu}
\bigl(\lambda^{-1}u(\lambda)-\mu^{-1}u(\mu)\bigr)\bigl(u(\lambda)-u(\mu)\bigr)
&=q\bigl(u(\lambda)-u(\mu)\bigr)\bigl(\lambda^{-1}u(\lambda)-
\mu^{-1}u(\mu)\bigr)\\
\label{relm}
\bigl(\lambda^{-1}m(\lambda)-\mu^{-1}m(\mu)\bigr)\bigl(m(\lambda)-m(\mu)\bigr)
&=q^{-1}\bigl(m(\lambda)-m(\mu)\bigr)\bigl(\lambda^{-1}m(\lambda)-
\mu^{-1}m(\mu)\bigr)\\
\label{relum}
u(\lambda)m(\mu)&=q^{-1}m(\mu)u(\lambda)
\end{align}
with $u(\lambda)=\sum_{i=0}^{\infty} (-1)^i u_{i+1}\lambda^{-i}$
and $m(\lambda)=\sum_{i=0}^{\infty} (-1)^i m_{i+1}\lambda^{-i}$
is a quantum deformation of the algebra of functions over the Poisson 
homogeneous
space $\bigl(H_-\backslash B_-,P_{\infty}\bigr)$.
\end{prop}
The following proposition shows that the action by vector fields
of the Abelian Lie subalgebra $\fh_+$ on $H_-\backslash B_-$ can 
also be quantized.
\begin{prop}\label{exih}
There is $(H_n)_{n\in\NN^*}$, a commutative family of derivations on
$\CC[H_-\backslash B_-]_q$ defined by the formulas~:
\begin{align}
\label{imul}
H_{\mu}\bigl(u(\lambda)\bigr)&=\dfrac{1}{\lambda^{-1}-\mu^{-1}}
\bigl(\lambda^{-1}u(\lambda)-\mu^{-1}u(\mu)\bigr)v(\mu)u(\lambda)
-\dfrac{\mu^{-1}}{\lambda^{-1}-\mu^{-1}}\bigl(u(\lambda)-u(\mu)\bigr)
\bigl(1+v(\mu)u(\mu)\bigr);\\
\label{imum}
H_{\mu}\bigl(m(\lambda)\bigr)&=\dfrac{\mu^{-1}}{\lambda^{-1}-\mu^{-1}}
\bigl(1+m(\mu)w(\mu)\bigr)
\bigl(m(\lambda)-m(\mu)\bigr)-\dfrac{1}{\lambda^{-1}-\mu^{-1}}m(\lambda)w(\mu)
\bigl(\lambda^{-1}m(\lambda)-\mu^{-1}m(\mu)\bigr)
\end{align}
with $H(\mu)=\sum_{k=1}^{\infty}(-1)^k H_k \mu^{-k}$,
$v(\mu)=-{\bigl(u(\mu)+\mu m(\mu)^{-1}\bigr)}^{-1}$ and 
$w(\mu)=-{\bigl(m(\mu)+\mu u(\mu)^{-1}\bigr)}^{-1}$.
This family of derivations deforms the classical action by vector fields of
$\fh_+$ on $H_-\backslash B_-$.
\end{prop}
The following theorem is a quantum version of the Drinfeld-Sokolov 
correspondence.
It shows that the quantization of the action of $\fh_+$ on 
$H_-\backslash B_-$
considered in Proposition \ref{exih} 
can be identified with the adjoint action 
of the integrals of motion on the phase space $A_q$.
\begin{thm}\label{thc}
There is an injective and birational map $\DS_q$ from
$\CC[H_-\backslash B_-]_q$
to $A_q$ defined by
\begin{equation}\label{fracu}
\DS_q(u(\lambda))=\lim\limits_{N\to\infty}
\dfrac{y_0^{-1}}{1+\dfrac{(\lambda x_0y_0)^{-1}}
      {1+\dfrac{(\lambda y_{-1}x_0)^{-1}}{\dfrac{\ddots
      (\lambda y_{-N+1}x_{-N+2})^{-1}}
      {1+(\lambda x_{-N+1}y_{-N+1})^{-1}}}}}
\end{equation}
and
\begin{equation}\label{fracm}
\DS_q(m(\lambda))=\lim\limits_{N\to\infty}
\dfrac{x_1^{-1}}{1+\dfrac{(\lambda x_1y_1)^{-1}}
      {1+\dfrac{(\lambda y_{1}x_2)^{-1}}{\dfrac{\ddots
      (\lambda y_{N-1}x_{N})^{-1}}
      {1+(\lambda x_{N}y_{N})^{-1}}}}}
\end{equation}
the limit being taken with respect to the $\lambda^{-1}$-adic topology, 
and where, 
in the first case, $\dfrac{a}{b}:=b^{-1}a$, and in the second case,
$\dfrac{a}{b}:=ab^{-1}$. We have~:
$\DS_q\circ H_n=\ad(I_n)\circ\DS_q$ for all integers $n$.
\end{thm}
The algebras $\CC[H_-\backslash B_-]_q$ and $A_q$ have fraction fields, and 
the phrase "birational" means that $DS_q$ induces an isomorphism of these
fraction fields.
\section{Commutativity of the local integrals 
of motion of the discrete sine-Gordon system}\label{unom}
We give here a (new) proof of the
commutativity of the quantum local integrals of motion 
\cite{IZKO} (see also \cite{FTT}, \cite{VOL}, \cite{HIK}).
This constitutes Proposition \ref{tha}.
Our proof is based on the explicit 
form taken by the elements $I_n$ of the basis of ${\cal I}$ given in the 
Proposition \ref{basi}.
First, we note that, as a consequence of Proposition \ref{nondema},
${\cal I}$ is a Lie subalgebra of ${\cal F}_0$.
We shall use Proposition \ref{expli} given below,
proved in \cite{Gru} and which is equivalent to Proposition \ref{basi}.
If $a$ and $b$ are two integers, we set~:
$$
\C{a}{b}=\left\{
          \begin{array}{ll}
           1,&\text{si }b=0\,;\\
           0,&\text{si }b<0\text{ ou si }b>\text{Max}(0,a)\,;\\
           \displaystyle{[a]!\over [b]! [a-b]!},&\text{si }0\leq b\leq a\,;
          \end{array}
         \right.
$$
with $[n]!=\prod_{i=1}^{n}[i]$ and $[i]=\dfrac{q^{i}-1}{q-1}$.\\
Also, for any relative integers $N,a_1,\ldots,a_N$, we set 
$F_q(a_1,\ldots,a_N)=\prod_{i=1}^{N}\C{a_i+a_{i+1}-1}{a_{i+1}}$.
\begin{prop}\label{expli}
A basis for ${\cal I}$ is given by the family $I_n=I(\psi_n)$, with 
$\psi_n=A_n+B_n$, $B_n=T^{\frac{1}{2}}A_n$, and
\begin{equation}\label{an}
A_n=\sum\frac{[n]}{[\alpha_1]}F_q(\alpha_1,\ldots,\alpha_{2N-2})
(x_1y_1)^{-\alpha_1}(y_1x_2)^{-\alpha_2}\ldots(y_{N-1}x_N)^{-\alpha_{2N-2}}
\end{equation}
the sum being taken on the set of indices $\a_1,\ldots,\a_{2N-2}$
such that $\a_i\in\NN$, $\a_1\in\NN^*$, and $\a_1+\ldots+\a_{2N-2}=n$.
Here, $N$ is any integer such that $n\leq 2(N-1)$.
\end{prop}
One will achieve the proof of the commutativity of ${\cal I}$ 
in several steps.
The key step is the existence of a filtration on a Lie subalgebra
${\cal F}'_0$ which contains ${\cal I}$.
\subsection{Gradation on ${\cal F}_0$}\label{defgra}
We denote $\deg_p$ the principal gradation on $A_q$
defined by $\deg_p(x_i)=\deg_p(y_i)=1$ for all $i\in\ZZ$.
Also, we set $e_{2i-1}=(x_iy_i)^{-1}$ and $e_{2i}=(y_ix_{i+1})^{-1}$.
The elements $e_j^{\pm 1}$, $j\in\ZZ$ generate $A_q[0]$. So, any element
$u\in{\cal F}_0$ can be represented by a sum of terms $P_k$, $k\in\ZZ$
with $\deg_p(P_k)=2k\in 2\ZZ$. Therefore, the gradation $\deg_p$
on $A_q$ induces a gradation $\deg$ on the module ${\cal F}_0$ by the following
way~: an element $u$ of ${\cal F}_0$ is said to be homogeneous of degree $k$ if there exists
$P\in A_q[0]$ such that $I(P)=u$ and $\deg_p(P)=2k$. The formula (\ref{alij1})
shows that $({\cal I},\deg)$ is a Lie subalgebra of ${\cal F}_0$
generated by $I_n=I(\psi_n)$, $n\in\NN^*$ with $\deg(I_n)=-n$.
\subsection{The subalgebra $B_q[0]$ of $A_q[0]$}
Let $B_q[0]$ be the subalgebra (without unit) of $A_q[0]$
generated by the $e_i$, $i\in\ZZ$. For all $i,j\in\ZZ$, we have 
$e_ie_{i+1}=q^{-1}e_{i+1}e_i$ and $e_ie_j=e_je_i$ if $|i-j|\geq 2$.
A basis for $B_q[0]$ is given by the family 
$\varepsilon_{\alpha}=\prod_{i\in\ZZ}e_i^{\alpha_i}$, where
$\a$ is a {\it{non-zero}} sequence in $\NN^{\ZZ}$ such that almost all
elements are zeros (except for a finite number of them).
We define a function $l$ on $B_q[0]$ by the following way.
If $P=\sum_{\alpha}\lambda_{\alpha}\varepsilon_{\alpha}$ is a non-zero element
in $B_q[0]$ we set $l(P):=\Inf\{ l(\alpha)/\,\lambda_{\alpha}\not= 0\}$
with $l(\alpha):=\card\{ i\in\ZZ/\, \alpha_i\not= 0\}$. We have
$l(P)\in\NN^*$. By convention, we consider that $l(0)=+\infty$.
Obviously, $B_q[0]$ is invariant by the translation automorphism
$T$ defined in \ref{defprin}, and $l\circ T=l$.
Moreover, if $P$ and $Q\in B_q[0]$, then
$l(PQ)\geq\Max\bigl( l(P),l(Q)\bigr)$ and
$l(P+Q)\geq\Inf\bigl( l(P), l(Q)\bigr)$ with equality 
in the last inequality if 
$l(P)\not= l(Q)$. 
\subsection{The Lie subalgebra ${\cal F}'_0$ of ${\cal F}_0$}\label{mar}
We note ${\cal F}'_0$ the quotient module $B_q[0]/\Im(T-Id)$.
A basis for ${\cal F}'_0$ is given by the elements 
$I(\varepsilon_{\alpha})$ where $\a$ is an almost zero sequence satisfying
the property $\alpha_i=0$ if $i\leq 0$,
and $\alpha_1\not= 0$ or $\alpha_2\not= 0$. The module ${\cal F}'_0$ is a Lie
subalgebra by taking the formula (\ref{alij1}) as a definition.
Furthermore, we have a natural injective map of Lie algebras
${\cal F}'_0\hookrightarrow{\cal F}_0$ such that the following 
diagram is commutative.
$$
\begin{array}{rcl}
B_q[0]&\hookrightarrow&A_q[0]\\
\downarrow&&\downarrow\\
{\cal F}'_0&\hookrightarrow&{\cal F}_0.
\end{array}
$$
The first horizontal map is the canonical injection of $B_q[0]$ in $A_q[0]$.
The vertical maps are the canonical projections. Note that
${\cal F}'_0$ is a graded Lie subalgebra of ${\cal F}_0$
with respect to gradation $\deg$. Also, according to Proposition \ref{expli},
${\cal F}'_0$ contains ${\cal I}$.
\subsection{Filtration on ${\cal F}'_0$}
Let us start with the following lemma.
\begin{lem}\label{lema}
Let $u\in{\cal F}'_0$ with $u\not=0$. Then 
$L_u:=\{ l(P)/\, P\in B_q[0]\, \text{and}\, I(P)=u\}$ 
is a bounded non-empty set in $\NN^*$.
\end{lem}

\begin{dem}
Obviously, $L_u\not=\emptyset$. Consider $Q\in B_q[0]$ such that $Q$ is generated by 
$e_i$, $n_1\leq i\leq n_2$.
Let $V_1$ (resp. $V_2$) be the submodule generated by the monomials $\varepsilon_{\a}$
with $l(\a)\leq n_2-n_1+1$ (resp. $l(\a)>n_2-n_1+1$). We have 
$B_q[0]=V_1\oplus V_2$, $Q\in V_1$, and $V_i$ ($i=1,2$) is invariant by $T$.
Assume that there exists $P\in B_q[0]$ such that $I(P)=u$ and $l(P)>n_2-n_1+1$.
Then, there exists $R$ such that $P=Q+T(R)-R$. We set $R=R_1+R_2$ with $R_i\in V_i$
($i=1,2$). By projecting on $V_2$, we obtain $P=T(R_2)-R_2$.
So, $u=0$. But this contradicts our hypothesis. 
\end{dem}

Lemma \ref{lema} allows us to define a length function on
${\cal F}'_0$.

\begin{Def}
{\it{ Let $u\in{\cal F}'_0$, with $u\not= 0$. We set $l(u)=\Max L_u$.
By convention, $l(0)=+\infty$.}} 
\end{Def}

Thanks to section \ref{mar} the following lemma allows us 
to compute lengths in ${\cal F}'_0$ explicitly.
\begin{lem}\label{lemb}
Let $x\in B_q[0]$ be a non-zero element such that $x$ is a linear combination
of monomials of the form $\varepsilon_{\a}$, with $\a_i= 0$ if $i\leq 0$
and $\a_1\not= 0$ or $\a_2\not= 0$.
Then, $l\bigl(I(x)\bigr)=l(x)$.
\end{lem}  

\begin{dem}
Set $l(x)=k$. Let $V_1$ (resp. $V_2$) denotes the submodule of $B_q[0]$
generated by all monomials $\varepsilon_{\a}$ with length $k$
(resp. with a length different from $k$). We proceed as in the proof of
the Lemma \ref{lema}
\end{dem}

We are now able to show that ${\cal F}'_0$ is a filtred Lie algebra.
\begin{lem}\label{lemc}
Let  $u$ and $v$ be two elements of ${\cal F}'_0$.
Then, $l\bigl([u,v]\bigr)\geq\Max\bigl(l(u),l(v)\bigr)$.
\end{lem}

\begin{dem}
Set $j=l(u)$, $k=l(v)$ and $n=\Max(j,k)$. There exist $P$ and $Q$
in $B_q[0]$ such that $I(P)=u$, $I(Q)=v$, $l(P)=j$ and $l(Q)=k$.
For all $k\in\ZZ$, we have $l(T^{\alpha}(P))=l(P)$. So,
$l\bigl((T^{\alpha}(P))Q\bigr)\geq n$ and $l\bigl(Q(T^{\alpha}(P))\bigr)\geq n$.
So, $l\bigl([T^{\alpha}(P),Q]\bigr)\geq n$.
Hence, $l\bigl(\sum_{\alpha=-\infty}^{\infty}[T^{\alpha}(P),Q]\bigr)\geq n$
and $l\bigl([u,v]\bigr)\geq n$.
\end{dem}

\subsection{End of the proof}
For $n\in\NN^*$, we set $u_n=e_1^n$ and $v_n=e_2^n$.
Thanks to Proposition \ref{expli},
there exists $w_n\in B_q[0]$ such that $l(w_n)\geq 2$ and 
$\psi_n=u_n+v_n+w_n$. Clearly we have
$[I(u_n),I(u_p)]=[I(v_n),I(v_p)]=0$ for all $n,p\in\NN$.
On the other hand, by using Lemma \ref{lemb},
the computation shows that $l\bigl([I(u_n),I(v_p)]\bigr)=2$.
So, we deduce from Lemma \ref{lemc} and the bilinearity of the Lie bracket that 
$l\bigl([I_n,I_p]\bigr)\geq 2$ for all $n,p\in\NN$.
But, for degree reasons, $[I_n,I_p]$ is proportional to $I_{n+p}$
and  $l(I_{n+p})=1$. Hence, $[I_n,I_p]=0$.

\section{Quantization of the Poisson homogeneous space 
$\bigl(H_-\backslash B_-,P_{\infty}\bigr)$}\label{secquanti}
The aim of this section is to prove Proposition \ref{thb}.
As shown in the subsection \ref{ahbon}, 
at the semi-classical level, the generators
$u_i$ and $m_i$, $i>0$, give a natural system of coordinate
functions on $H_-\backslash B_-$ which satisfy the Poisson relations
obtained by taking the limit $q\to 1$ in 
(\ref{relu}), (\ref{relm}), (\ref{relum}).
Therefore, to prove Proposition \ref{thb}, it is enough to show that
$\CC[H_-\backslash B_-]_q$ is a flat deformation of $\CC[H_-\backslash B_-]$.
The idea is to obtain a realization of the algebra
$\CC[H_-\backslash B_-]_q$ in $A_q$ by using Lemma \ref{enrf}
which asserts that 
the finite screening charges of the discrete sine-Gordon system satisfy
the quantum Serre relations. In all the following, the ground ring
is no longer $\QQ[q,q^{-1}]$ but $\CC[q,q^{-1}]$.

\subsection{The Poisson homogeneous space 
$\bigl(H_-\backslash B_-,P_{\infty}\bigr)$}\label{ahbon}
The Poisson manifold $\bigl(H_-\backslash B_-,P_{\infty}\bigr)$
was defined in \ref{defprin}. Any element 
${\bar x}\in H_-\backslash B_-$ can be expressed uniquely in the 
following form~: 
$$
{\bar x}=\cl_{H_-}
\begin{pmatrix}
1&v_{\cl}(\lambda)({\bar x})\\
0&1
\end{pmatrix}
\begin{pmatrix}
1&0\\
u_{\cl}(\lambda)({\bar x})&1
\end{pmatrix}
$$
with $u_{\cl}(\lambda)\in\CC[H_-\backslash B_-][[\lambda^{-1}]]$
and $v_{\cl}(\lambda)\in\lambda^{-1}\CC[H_-\backslash B_-]
[[\lambda^{-1}]]$.
For $i\in\NN^{*}$, we define  coordinate functions $u_{i,\cl}$
and $m_{i,\cl}$ on $H_-\backslash B_-$
by~:
\begin{align}
u_{\cl}(\lambda)&=
\sum_{i=0}^{\infty}
(-1)^i u_{i+1,\cl}\lambda^{-i}\\
m_{\cl}(\lambda)&=
\sum_{i=0}^{\infty}
(-1)^i m_{i+1,\cl}\lambda^{-i}\\
\textnormal{and}\quad m_{\cl}(\lambda)&:=
-\lambda v_{\cl}(\lambda)
\bigl(1+u_{\cl}
v_{\cl}(\lambda)\bigr)^{-1}\in\CC[H_-\backslash B_-][[\lambda^{-1}]].
\end{align}
The functions $u_{i,\cl}$ and $m_{i,\cl}$ are 
algebraically independent and 
$\CC[H_-\backslash B_-]=\CC[u_{i,\cl},m_{i,\cl},i>0]$.
Moreover, computation shows that the Poisson relations between these
functions (the Poisson structure is induced by the field
of bivectors $P_{\infty}$) are precisely the ones we get from
(\ref{relu}), (\ref{relm}), (\ref{relum})
when $q\to 1$.

\subsection{The Enriquez-Feigin morphism}
Let $\fn_-$ be a nilpotent subalgebra of ${\widehat{\mathfrak{sl}_2}}$
and $U_q \fn_-$ be the quantum algebra  generated by the generators 
$f_+$ and $f_-$ subject to the quantum Serre relations~:
\begin{equation}\label{serreq}
f_{\pm}^3 f_{\mp}-(q+1+q^{-1})(
f_{\pm}^2 f_{\mp}f_{\pm}-f_{\pm}f_{\mp}f_{\pm}^2)
-f_{\mp}f_{\pm}^3=0.
\end{equation}
Let $\deg$ be the gradation on $U_q \fn_-$ defined by 
$\deg f_{\pm}=\pm 1$. In the sequel of the article,
if $(A,\deg_A)$ and $(B,\deg_B)$ are two graded algebras, we define
the twisted tensor product $A\ootimes B$ by the formula~:
$$
(a_1\ootimes b_1)(a_2\ootimes b_2)=
q^{-\deg_A(a_2)\deg_B(b_1)}(a_1 a_2\ootimes b_1 b_2)
$$
for homogeneous elements $a_1,a_2,b_1,b_2$ in $A$ and $B$.
There is a unique graded algebra morphism $\DDelta$ from
$U_q\fn_-$ to $U_q\fn_-\ootimes U_q\fn_-$ given by
$\DDelta(f_{\pm})=f_{\pm}\ootimes 1+1\ootimes f_{\pm}$.
This morphism is called the twisted comultiplication on 
$U_q\fn_-$.

\begin{lem}[\cite{ENR}]\label{enrf}
Let $n\in\NN$. Then $\SPN:=\sum_{i=1}^{n}x_i$
and $\SMN:=\sum_{i=1}^{n}y_i$ satisfy the quantum Serre relations.
In other words, there exists a graded algebra morphism $f_n$
defined by
\begin{equation}\label{morphie}
\begin{array}{rcl}
f_n~:\quad U_q\fn_-&\longrightarrow&A_q\\
f_{\pm}&\longmapsto&\SPMN.
\end{array}
\end{equation}
\end{lem}

\begin{dem}
For $1\leq i\leq n$, we define two graded algebra morphisms
$\varphi_i$ and $\psi_i$ by~:
$$
\begin{array}{rcl}
\varphi_i~:\quad U_q\fn_-&\longrightarrow&\CC[x_i,x_i^{-1}]\\
f_+&\longmapsto&x_i\\
f_-&\longmapsto&0
\end{array}
$$
and 
$$
\begin{array}{rcl}
\psi_i~:\quad U_q\fn_-&\longrightarrow&\CC[y_i,y_i^{-1}]\\
f_+&\longmapsto&0\\
f_-&\longmapsto&y_i
\end{array}
$$
with the convention $\deg x_i=-\deg y_i =1$.
Then, it appears that the map $f_n$ defined in (\ref{morphie})
is equal to 
$\mult_{2n}\circ (\varphi_1\ootimes\psi_1\ootimes\ldots\ootimes
\varphi_n\ootimes\psi_n)\circ \DDelta^{(2n)}$
where $\mult_{2n}$ is the algebras monomorphism
$$
\begin{array}{rcl}
\mult_{2n}~:\quad \CC[x_1,x_1^{-1}]\ootimes
\CC[y_1,y_1^{-1}]\ootimes
\CC[x_n,x_n^{-1}]\ootimes
\CC[y_n,y_n^{-1}]&\longrightarrow&A_q\\
u_1\ootimes v_1\ootimes \ldots u_n\ootimes v_n&\longmapsto&
u_1v_1\ldots u_n v_n.
\end{array}
$$
\end{dem}

So as to realize $\CC[H_-\backslash B_-]_q$
in $A_q$, it is useful to give another expression of the quantum
algebra $U_q\fn_-$.

\begin{Def}
{\it{Let $q^{\frac{1}{4}}$ be an indeterminate, 
$q^{\frac{1}{2}}=\bigl(q^{\frac{1}{4}}\bigr)^{2},\,
 q=\bigl(q^{\frac{1}{2}}\bigr)^{2},\,
 K_0=\CC[q^{\frac{1}{4}},q^{-\frac{1}{4}}]$, $K=\CC[q,q^{-1}]$,
 \begin{align}
  \label{mat2}
   H&= \begin{pmatrix}
         q^{-\frac{1}{}}&0&0&0\\
         0&q^{\frac{1}{4}}&0&0\\
         0&0&q^{\frac{1}{4}}&0\\
         0&0&0&q^{-\frac{1}{4}}
       \end{pmatrix} \in M_4(K_0)\simeq M_2(K_0)^{{\otimes}^2},\\
 \label{mat1}
  \text{and } R({\lambda,\mu})&=
   \begin{pmatrix}
     1&0&0&0\\
     0&\frac{\lambda -\mu}{q^{-\frac{1}{2}}\lambda -q^{\frac{1}{2}}\mu}&
     \frac{(q^{-\frac{1}{2}}-q^{\frac{1}{2}})\mu}{ 
     q^{-\frac{1}{2}}\lambda -q^{\frac{1}{2}}\mu}&0\\
     0&\frac{(q^{-\frac{1}{2}}-q^{\frac{1}{2}})\lambda}{ 
     q^{-\frac{1}{2}}\lambda -q^{\frac{1}{2}}\mu}&
     \frac{\lambda -\mu}{q^{-\frac{1}{2}}\lambda -q^{\frac{1}{2}}\mu}&0\\
     0&0&0&1
   \end{pmatrix} \in M_4(K_0)\simeq M_2(K_0)^{{\otimes}^2}.
 \end{align}
By definition, $\CC[N_+]_q$ is the algebra generated over $K$
by the $a_{i,j}^{(r)}$ for $i,j\in\{1,2\}\, r\in\NN$
and the relations~:
 \begin{gather}
 \label{RR1}
 a_{2,1}^{(0)}=0,
 a_{1,1}^{(0)}=
 a_{2,2}^{(0)}=1,\\
 \label{RR2}
 R({\lambda,\mu}){\cal L}^1(\lambda)H{\cal L}^2(\mu)=
 {\cal L}^2(\mu)H{\cal L}^1(\lambda)R({\lambda,\mu}),\\
 \label{RR3}
 a_{1,1}(q\lambda)\bigl[
 a_{2,2}(\lambda)-
 a_{2,1}(\lambda)a_{1,1}(\lambda)^{-1}a_{1,2}(\lambda)\bigr]=1,
 \end{gather}
 with~:
 $a_{i,j}(\lambda)=\sum\limits_{r=0}^{+\infty}a_{i,j}^{(r)}
 \lambda^r,\,
 {\cal L}(\lambda)=[a_{i,j}(\lambda)],
 {\cal L}^1(\lambda)={\cal L}(\lambda)\otimes Id$
 and 
 ${\cal L}^2(\mu)=Id\otimes {\cal L}(\mu)$.}}
\end{Def}

Relations (\ref{RR2}) and (\ref{RR3})
have coefficients in $K$. So, the above definition makes sense.
The relation (\ref{RR3}) is the quantum determinant relation. 
It can be shown \cite{Grut} that $\CC[N_+]_q$ is a quantum deformation
of the algebra of functions on the Poisson manifold $N_+$ defined in 
\ref{defprin}, equipped with the Poisson bivector
$P=r^G-r^D+\dfrac{1}{4}(h^G\otimes h^D-h^D\otimes h^G)$
where $r$ denotes the $r$-trigonometric matrix~:
$$
r(\lambda,\mu)=-\dfrac{1}{4}\dfrac{\lambda+\mu}{\lambda-\mu}h\otimes h
 -\dfrac{1}{\lambda-\mu}(\lambda e\otimes f+\mu f\otimes e)
$$
with $h=\begin{pmatrix}
1&0\\
0&-1\end{pmatrix}$,
$e=\begin{pmatrix}
0&1\\
0&0
\end{pmatrix}$ and $f=\begin{pmatrix}
0&0\\
1&0\end{pmatrix}$.
We can define a gradation $\deg$ on $\CC[N_+]_q$
by $\deg(a_{i,j}^{(k)})=i-j$ for integers $i,j,k$
as well as a twisted comultiplication $\DDelta$ given by
$\DDelta\bigl( {\cal L}(\lambda)\bigr)=
{\cal L}(\lambda)\ootimes {\cal L}(\lambda)$.
Moreover, the map 
$$
 \begin{array}{rcl}
 \Phi~:\quad U_q\nm&\longrightarrow&\CC[N_+]_q\\
 f_{+}&\longmapsto&a_{2,1}^{(1)}\\
 f_{-}&\longmapsto&a_{1,2}^{(0)}
 \end{array}
$$
is an algebra morphism such that 
$\DDelta\circ\Phi=(\Phi\ootimes\Phi)\circ\DDelta$.
Let us consider the morphism $f_n$ defined 
in Lemma \ref{enrf} above.
Thanks to the proof of Lemma \ref{enrf}, 
there exists a graded algebra morphism
\begin{equation}\label{petit}
 \begin{array}{rcl}
 g_n~:\quad \CC[N_+]_q&\longrightarrow&A_q\\
 {\cal L}(\lambda)&\longmapsto&\prod_{i=1}^{n}\begin{pmatrix}
 1&0\\
 \lambda x_i&1\end{pmatrix}\begin{pmatrix}
 1&y_i\\
 0&1\end{pmatrix}.
 \end{array} 
\end{equation}
It is clear that for $r>n$, 
$a_{1,1}^{(r-1)},\, a_{1,2}^{(r-1)},\, a_{2,1}^{(r)},\, 
a_{2,2}^{(r)}\in\ker g_n$. Moreover,
$a_{1,1}^{(n-1)},\, a_{1,2}^{(n-1)},\, a_{2,1}^{(n)}$
and $a_{2,2}^{(n)}$ are invertible.
This leads to the study of the quantum algebra
$\CC[B_- w_n B_-\cap N_+]_q$ given below.

\subsection{The quantum Schubert cell 
$\CC[B_- w_n B_-\cap N_+]_q$}\label{cell}
The interest in studying the algebra $\CC[B_- w_n B_-\cap N_+]_q$
stems from the fact that the generating series of certain functions
defined on this quantum algebra satisfies the same relation 
(\ref{relu}) as the generators $u_i$, $i\in\NN^*$ in 
$\CC[H_-\backslash B_-]_q$, and that we can deduce from 
(\ref{petit}), the existence of an algebra morphism
from $\CC[B_- w_n B_-\cap N_+]_q$ to $A_q$.

\begin{Def}
{\it{The algebra $\CC[B_- w_n B_-\cap N_+]_q$ 
is given by generators 
$a_{i,j}^{(k)}\, (i,j\in\{ 1,2\}; k\in\{ 0,\ldots,n\}),\, {a_{1,1}^{(n-1)}}',\,
{a_{1,2}^{(n-1)}}',\, {a_{2,1}^{(n)}}',\,
{a_{2,2}^{(n)}}'$,
and relations (\ref{RR1}), (\ref{RR2}), (\ref{RR3})
with $a_{i,j}(\lambda)=\sum_{k=0}^{n}a_{i,j}^{(k)}\lambda^k$
and $a_{1,1}^{(n)}=a_{1,2}^{(n)}=0$ as well as
relations which express the fact that ${a_{1,1}^{(n-1)}}'$
(resp. ${a_{1,2}^{(n-1)}}'$,\, ${a_{2,1}^{(n)}}'$,\, ${a_{2,2}^{(n)}}'$)
is an inverse for $a_{1,1}^{(n-1)}$
(resp. $a_{1,2}^{(n-1)}$,\, $a_{2,1}^{(n)}$,\, $a_{2,2}^{(n)}$).}}
\end{Def}

At the semi-classical limit, for $q\to 1$, we get the algebra
of functions on the Schubert cell 
$(B_- w_n B_-\cap N_+,P)$ with $w_n=\diag(\lambda^{-n},\lambda^n)$
with a Poisson structure given 
by the fact that $(B_- w_n B_-\cap N_+,P)$ may be viewed
as a symplectic leaf of the Poisson manifold
$(N_+,P)$. 
The algebra $\CC[B_- w_n B_-\cap N_+]_q$ is just a rough quantum deformation
of $(B_- w_n B_-\cap N_+,P)$. To obtain an exact quantum deformation,
one should impose relations between 
${a_{1,1}^{(n-1)}}',\,
{a_{1,2}^{(n-1)}}',\, {a_{2,1}^{(n)}}',\,
{a_{2,2}^{(n)}}'$ and all the $a_{i,j}^{(k)}$ on the definition.
However, we don't need to be so precise and our definition will suffice.
There is a natural morphism 
$p:\CC[N_+]_q\longrightarrow\CC[B_- w_n B_-\cap N_+]_q$.
Moreover, if $C$ is an algebra and if 
$f:\CC[N_+]_q\longrightarrow C$ is an algebra morphism
such that 
$a_{1,1}^{(r-1)},\, a_{1,2}^{(r-1)},\,
a_{2,1}^{(r)},\, a_{2,2}^{(r)}\in\ker f$ pour $r>n$,
and $f\bigl(a_{1,1}^{(n-1)}\bigr),\, f\bigl(a_{1,2}^{(n-1)}\bigr),\,
f\bigl(a_{2,1}^{(n)}\bigr),\, f\bigl(a_{2,2}^{(n)}\bigr)$ 
are invertible in
$C$, then there exists $g~:\CC[B_- w_n B_-\cap N_+]_q\longrightarrow C$
an algebra morphism such that the following diagram is
commutative~:
$$
 \begin{array}{rcl}
 \CC[N_+]_q&{\overset{p}{\longrightarrow}}&\CC[B_- w_n B_-\cap N_+]_q\\
 f\searrow&&\swarrow g\\
 &C&
 \end{array}
$$
By virtue of (\ref{petit}), it follows that 
there exists an algebra morphism $h_n:\CC[B_- w_n B_-\cap N_+]_q\longrightarrow A_q$
such that~:
\begin{equation}\label{hnl}
h_n\bigl(L(\lambda)\bigr)=\prod_{i=1}^{n}\begin{pmatrix}
 1&0\\
 \lambda x_i&1\end{pmatrix}\begin{pmatrix}
 1&y_i\\
 0&1\end{pmatrix}
\end{equation}
with $L(\lambda)=[a_{i,j}(\lambda)]_{i,j\in\{ 1,2\}}$
and $a_{i,j}(\lambda)=\sum_{k=0}^{n}a_{i,j}^{(k)}\lambda^{k}$.
The element
$a_{2,2}(\lambda)$ is invertible in the ring 
$\CC[B_- w_n B_-\cap N_+]_q((\lambda^{-1}))$. 
We set $\a(\lambda):=a_{2,2}(\lambda)^{-1}a_{2,1}(\lambda)$.
\begin{lem}\label{relal}
The function $\alpha(\lambda)$ satisfies the same relation
(\ref{relu}) as the function $u(\lambda)$.
We have~:
$\bigl(\lambda^{-1}\alpha(\lambda)-
\mu^{-1}\alpha(\mu)\bigr)\bigl(\alpha(\lambda)-\alpha(\mu)\bigr)
=q\bigl(\alpha(\lambda)-\alpha(\mu)\bigr)\bigl(\lambda^{-1}\alpha(\lambda)-
\mu^{-1}\alpha(\mu)\bigr)$.
\end{lem}

\begin{dem}
By definition, for two elements $a$ and $b$ of an algebra, 
We denote by [a,b] the commutator ab - ba.
Then,
\begin{align*}
\bigl[\alpha(\lambda),\alpha(\mu)\bigr]&=
{a_{2,2}(\lambda)}^{-1}
\bigl[ a_{2,1}(\lambda),{a_{2,2}(\mu)}^{-1}\bigr] 
a_{2,1}(\mu)+
\bigl[{a_{2,2}(\lambda)}^{-1},{a_{2,2}(\mu)}^{-1}\bigr]\,
a_{2,1}(\lambda)a_{2,1}(\mu)\\
&+{a_{2,2}(\mu)}^{-1}{a_{2,2}(\lambda)}^{-1}
\bigl[ a_{2,1}(\lambda),a_{2,1}(\mu)\bigr]
+{a_{2,2}(\mu)}^{-1}
\bigl[
{a_{2,2}(\lambda)}^{-1},a_{2,1}(\mu)
\bigr]\,
a_{2,1}(\lambda).
\end{align*}
Relation (\ref{RR2}) shows that 
\begin{align}
\label{taga}
\forall\, i,j\in\{ 1,2\},\quad [a_{i,j}(\lambda),a_{i,j}(\mu)]&=0\\
\label{tagb}
\text{and}\quad \bigl[a_{2,1}(\lambda),a_{2,2}(\mu)\bigr]&=
\bigl[a_{2,1}(\mu),a_{2,2}(\lambda)\bigr]\\
\label{tagc}
&=(1-q^{-1})\bigl(\mu a_{2,2}(\mu)a_{2,1}(\lambda)
-\lambda a_{2,2}(\lambda)a_{2,1}(\mu)\bigr).
\end{align}
So, thanks to (\ref{taga}) and (\ref{tagb}), we get
\begin{align}
\bigl[\alpha(\lambda),\alpha(\mu)\bigr]&=a_{2,2}(\lambda)^{-1}
a_{2,2}(\mu)^{-1}\bigl[a_{2,1}(\lambda),a_{2,2}(\mu)\bigr]
\bigl(a_{2,2}(\lambda)^{-1}a_{2,1}(\lambda)
-a_{2,2}(\mu)^{-1}a_{2,1}(\mu)\bigr)\\
&=(1-q^{-1})\bigl(\mu \alpha(\lambda)-\lambda\alpha(\mu)\bigr)
\bigl(\alpha(\lambda)-\alpha(\mu)\bigr).
\end{align}
The result follows from this last equality.
\end{dem}

Let us see now what is the image
of $\a(\lambda)$ by the map $h_n$.

\begin{lem}\label{ima}
We have $h_n\bigl(\alpha(\lambda)\bigr)=
\dfrac{y_n^{-1}}{1+\dfrac{(\lambda x_n y_n)^{-1}}
       {1+\dfrac{(\lambda y_{n-1}x_n)^{-1}}{\dfrac{\ddots
       (\lambda y_{1}x_{2})^{-1}}
       {1+(\lambda x_{1}y_{1})^{-1}}}}}$,
with $\dfrac{a}{b}:=b^{-1}a$ for two elements $a$ and $b$ 
such that $b$ is invertible.
\end{lem}

\begin{dem}
Let us define elements $a_n, b_n, c_n, d_n$ in $A_q[\lambda]$
such that 
$h_n\bigl(L(\lambda)\bigr)=
  \begin{pmatrix}
    a_n&b_n\\
    c_n&d_n
  \end{pmatrix}$.
Then, thanks to (\ref{hnl}),
we have 
$h_n\bigl(a_{2,1}(\lambda)\bigr)=c_n$,
$h_n\bigl(a_{2,2}(\lambda)\bigr)=d_n$
and $h_n\bigl(\alpha(\lambda)\bigr)=d_n^{-1}c_n$. 
If $n=1$ then $d_1=1+\lambda(x_1y_1)$ is invertible 
in $A_q((\lambda^{-1}))$ and
$h_1\bigl(\alpha(\lambda)\bigr)=
\bigl(1+(\lambda x_1y_1)^{-1}\bigr)^{-1}y_1^{-1}$.
Let us make the hypothesis that $n>1$.
We have~:
\begin{align}
\label{kak}
c_n&=c_{n-1}+\lambda d_{n-1} x_n\\
\label{kok}
d_n&=c_{n-1}y_n +d_{n-1}(\lambda x_n y_n +1).
\end{align}
{}From this, the computation shows that
$$
h_n\bigl(\alpha(\lambda)\bigr)
=\Bigl ( 1+ (\lambda x_n y_n)^{-1}\bigl ( 
1+q^{-1}(d_{n-1}^{-1}c_{n-1}y_{n-1})(\lambda y_{n-1}x_n)^{-1}
\bigr )^{-1}\Bigr )^{-1}y_n^{-1}
$$
Therefore, if we set $V_n:=h_n\bigl(\alpha(\lambda)\bigr) y_n$ and
$W_n^{(k)}:=\bigl ( 1+q^{k} V_{n-1}(\lambda y_{n-1} x_n)^{-1}
\bigr )^{-1}$,
we see that 
$$
V_n=\bigl ( 1+ (\lambda x_n y_n)^{-1}W_n^{(-1)}\bigr )^{-1}.
$$
By induction on $n$, we deduce that 
for all integers $k$ with $k>n$,
$(x_k y_k)$ and $V_n$ commute.
Hence, by induction on $p$, 
$$
\forall p\in{\mathbb N},\quad
(\lambda x_n y_n)^{-p}W_n^{(k)}=W_n^{(k+p)}(\lambda x_n y_n)^{-p}.
$$
So,
\begin{align*} 
V_n&=\bigl ( 1+W_n^{(0)}(\lambda x_n y_n)^{-1}\bigr )^{-1}\\
&=\Bigl (1+\bigl ( 1+V_{n-1}(\lambda y_{n-1} x_n)^{-1}\bigr )^{-1}
(\lambda x_n y_n)^{-1}\Bigr )^{-1}.
\end{align*}
The result follows from this by induction on $n$.
\end{dem}

Set $\alpha(\lambda)=\sum_{i=0}^{\infty}(-1)^i\alpha_{i+1}\lambda^{-i}$
and let us see what the images of $\a_i$ in $A_q$ are.
For that, we shall need the following proposition.
\begin{prop}\label{aba}
Let $N$ be an integer and $I_N$ the ideal in 
$A_N:=\CC[q,q^{-1}]\{\{ t_1,\ldots,t_N\}\}$
generated by elements $t_i t_{i+1}-q t_{i+1}t_i$ for
$i\in\{ 1,\ldots,N\}$ and 
$t_it_j-t_j t_i$ for $|i-j|\geq 2$.
Then, in the ring $A_N/I_N$, we have~: 
    $$
       \biggl (1-\Bigl (1-\bigl (1-\ldots (1-t_N)^{-1}t_{N-1}\bigr )^{-1}
       \ldots t_2\Bigr )^{-1}t_1 \biggr )^{-1}
       =\displaystyle\sum_{\alpha_1,\ldots,\alpha_N}
       F_q(\alpha_1,\ldots,\alpha_N)t_N^{\alpha_N}\ldots t_1^{\alpha_1}.
    $$
We recall that the function $F_q$ has been defined in section 
\ref{unom}.
\end{prop}

\begin{dem}
We note ${\cal F}_N$ the quantum fraction in the left hand side,
$v_N$ the valuation on $A_N$ corresponding to the gradation given by
$\deg t_j$ for all $j$, and $i_N$ the valued injection from
$A_N$ to $A_{N+1}$ given by $i_N(t_j)=t_{j+1}$. 
If $N=1$, the result is obvious.
Let us assume that the property is true until the rank $N$. Then,
$v_{N+1}\bigl(i_N({\cal F}_N)t_1\bigr)\geq 1$ for $v_N({\cal F}_N)\geq 0$.
So, $1-i_N({\cal F}_N)t_1$ is invertible 
in $A_{N+1}$ and ${\cal F}_{N+1}$ exists. 
Set $t'_1=t_1,\ldots,t'_{N-1}=t_{N-1}$
and $t'_N=(1-t_{N+1})^{-1}t_N$.
Then, the $t'_j$ satisfy the same relations as the generators
$t_j$ in $A_N/I_N$. Moreover,
$$
{\cal F}_{N+1}=\biggl 
(1-\Bigl (1-\bigl (1-\ldots (1-{t'}_N)^{-1}{t'}_{N-1}\bigr )^{-1}
\ldots {t'}_2\Bigr )^{-1}{t'}_1 \biggr )^{-1}.
$$
So, the induction hypothesis implies that
${\cal F}_{N+1}=\displaystyle\sum_{\alpha_1,\ldots,\alpha_N}
F_q(\alpha_1,\ldots,\alpha_N)\, {t'}_N^{\alpha_N}\ldots {t'}_1^{\alpha_1}$.
On the other hand, an induction on $k$ shows that~:
$$
\bigl [(1-t_{N+1})^{-1}t_N\bigr ]^{k}=
(1-t_{N+1})^{-1}\ldots (1-q^{k-1}t_{N+1})^{-1}t_N^{k}.
$$
Therefore,
$$
{\cal F}_{N+1}=\displaystyle\sum_{\alpha_1,\ldots,\alpha_N}
F_q(\alpha_1,\ldots,\alpha_N)
(1-t_{N+1})^{-1}\ldots (1-q^{\alpha_{N}-1}t_{N+1})^{-1}
t_{N}^{\alpha_N}\ldots {t}_1^{\alpha_1}.
$$
The result follows from the classical relation~:
$$
\displaystyle\prod_{s=0}^{N-1}(1-q^st)^{-1}
=\displaystyle\sum_{k\geq 0}\C{N+k-1}{k}t^{k}.
$$
\end{dem}

Lemma \ref{ima} and Proposition \ref{aba} allow us to 
obtain explicitly images of
$\a_i$ by $h_n$.
\begin{cor}\label{uncor}
We have 
$h_n(\alpha_i)=
\sum F_q(\alpha_{1},\ldots,\alpha_{2n-1})(x_{1}y_{1})^{-\alpha_{2n-1}}\ldots
(x_n y_n)^{-\alpha_{1}}y_n^{-1}$,
the sum being taken on all integers $\alpha_1,\ldots,\alpha_{2n-1}$
such that $\alpha_1+\ldots+\alpha_{2n-1}=i-1$.
\end{cor}

The fact that the $\a_i$ satisfy relation 
(\ref{relu}) leads to the study of the following quantum algebra.
\subsection{The quantum homogeneous space
$\CC[S_{\infty}\backslash B_-]_q$}
Let $S_{\infty}$ be the sub-group of $B_-$ constituted by all lower
triangular matrices of the form 
$\begin{pmatrix}
 a&\lambda^{-1}b\\
 0&a^{-1}
\end{pmatrix}.$

\begin{Def}
{\it{We denote by $\CC[S_{\infty}\backslash B_-]_q$ the algebra
given by generators $u_i,\, i\in\NN^{*}$, 
and relation (\ref{relu})~:
$
\bigl(\lambda^{-1}u(\lambda)-\mu^{-1}u(\mu)\bigr)\bigl(u(\lambda)-u(\mu)\bigr)
=q\bigl(u(\lambda)-
u(\mu)\bigr)\bigl(\lambda^{-1}u(\lambda)-\mu^{-1}u(\mu)\bigr)$
with $u(\lambda)=\sum_{i=0}^{\infty}(-1)^i u_{i+1}\lambda^{-i}$.}}
\end{Def}

We can check that relations coming from (\ref{relu})
are equivalent to the equalities
\begin{equation}\label{uij}
\forall\, i<j,\quad [u_i,u_j]=(1-q^{-1})\sum_{k=i}^{i+j-1}u_k u_{i+j-k}.
\end{equation}
The algebra $\CC[S_{\infty}\backslash B_-]_q$ is a graded algebra
with the gradation given by $\deg u_i=1$ for all $i$.
Thanks to Lemma \ref{relal}, there exists a specialization
morphism~:
\begin{equation}
\begin{array}{rcl}
r~:\quad\CC[S_{\infty}\backslash B_-]_q&\longrightarrow&
\CC[B_- w_n B_-\cap N_+]_+\\
\forall\, i,\quad u_i&\longmapsto&\alpha_i
\end{array}.
\end{equation}
We set $h'_n=T^{-n}\circ h_n\circ r$ where
$T$ is the translation automorphism on $A_q$.
Thanks to Corollary ~\ref{uncor},
for all integers $i,j,m$ with $i\leq 2n$ and $i\leq 2m$
we have $h'_n(u_i)=h'_m(u_i)$. 
Now, if we take into account Lemma \ref{ima},
we deduce the existence of a graded algebra morphism
\begin{equation}\label{morh}
\begin{array}{rcl}
h~:\quad \CC[S_{\infty}\backslash B_-]_q&\longrightarrow& A_q\\
u(\lambda)&\longmapsto&\lim\limits_{N\to\infty}
\dfrac{y_0^{-1}}{1+\dfrac{(\lambda x_0y_0)^{-1}}
       {1+\dfrac{(\lambda y_{-1}x_0)^{-1}}{\dfrac{\ddots
       (\lambda y_{-N+1}x_{-N+2})^{-1}}
       {1+(\lambda x_{-N+1}y_{-N+1})^{-1}}}}}
 \end{array}
\end{equation}
with the convention that $\dfrac{a}{b}=b^{-1}a$ if $b$ is invertible.
Explicitly, the image of $u_i$ by $h$ is given by the formula~:
\begin{equation}\label{foru}
h(u_i)=\sum F_q(a_1,a_2,\ldots)\ldots
(x_{-k}y_{-k})^{-a_{2k+1}}(y_{-k}x_{-k+1})^{-a_{2k}}
\ldots (x_0 y_0)^{-a_1}y_0^{-1},
\end{equation}
the sum being taken  on all integers $a_i$ such that
$\sum_k a_k=i-1$. For example, $h(u_1)=y_0^{-1}$ and
$h(u_2)=(x_0 y_0)^{-1}y_0^{-1}$.
We note that for all integers $i$ and $n$ with $i\leq 2(n-1)$
we have $h_n(\a_i)=T^{n}h(u_i)$.

\begin{prop}\label{deff}
A basis for $\CC[S_{\infty}\backslash B_-]_q$ is given
by the family
$\xi_a:=\prod_{i=1}^{\infty}u_i^{a_i}$ where
$a=(a_i)_{i\in\NN^*}$
denotes an almost zero sequence of integers.
\end{prop}

\begin{dem}
Thanks to (\ref{uij}), $\CC[S_{\infty}\backslash B_-]_q$ 
is spanned by the family $\xi_a$. 
But this set of vectors is also free.
Indeed, this is a consequence of (1) the existence of $h$ given above,
(2) the fact that in each new element of the sequence $h(u_i)$ 
occurs one and only new element of the form $x_{-k}$ or
$y_{-k}$ (according to the parity of $i$) and (3)
the fact that the family  
$\prod x_i^{a_i} y_i^{b_i}$ forms a basis of $A_q$
where $(a_i)$ and $(b_i)$ are almost zero sequences in 
$\ZZ^{\ZZ}$. 
\end{dem}

Note that the proof of Proposition \ref{deff} shows
the following result.

\begin{cor}\label{hinj}
The algebra morphism $h$ is injective.
\end{cor}

In the classical case, when $q\to 1$, we see from (\ref{foru}),
that $h$ is a birational map. 
We can also deduce from Proposition \ref{deff} that
$\CC[S_{\infty}\backslash B_-]_q$ is a quantum deformation
of the algebra of functions on the Poisson manifold
$S_{\infty}\backslash B_-$ equipped with the Poisson structure induced
by the field of bivectors $r^L-r^R$  \cite{Grut}. 

\subsection{End of the proof}
It is based on Proposition \ref{deff} and Lemma 
\ref{doub} which show together that
$\CC[H_-\backslash B_-]_q$ is in a way the quantum ``double''
of $\CC[S_{\infty}\backslash B_-]_q$.

\begin{Def}
{\it{We note 
$\CC[S_{\infty}\backslash B_-]_q^+$
the algebra given  by generators~:
$m_i,\, i\in\NN^*$
and relations~: (\ref{relm}), i.e.,
$
\bigl(\lambda^{-1}m(\lambda)-\mu^{-1}m(\mu)\bigr)\bigl(m(\lambda)-m(\mu)\bigr)
=q^{-1}\bigl(m(\lambda)-m(\mu)\bigr)\bigl(\lambda^{-1}m(\lambda)-\mu^{-1}m(\mu)\bigr)
$
with 
$m(\lambda)=\sum_{i=0}^{\infty}(-1)^i m_i\lambda^{-i}$.}}
\end{Def}
We note also by $\deg$ the gradation on $\CC[S_{\infty}\backslash B_-]_q^+$
defined by $\deg m_i=-1$ for all $i$ and by $\varphi^+$
the anti-isomorphism of algebras~:
\begin{equation}\label{deffi+}
\begin{array}{rcl}
\varphi^+~:\quad\CC[S_{\infty}\backslash B_-]_q&
\longrightarrow&\CC[S_{\infty}\backslash B_-]_q^+\\
\forall\, i\in\NN^*,\quad u_i&\longmapsto&m_i
\end{array}.
\end{equation}
The map $\varphi^+$ is an anti-graded involution.
On the other hand, there exists also an involution $\varphi$
on $A_q$ which is an
algebras anti-graded anti-automorphism defined by~:
\begin{equation}\label{definfi}
\begin{array}{rcl}
\varphi~:\quad A_q&\longrightarrow&A_q\\
\forall\, i\in\ZZ,\quad x_i&\longmapsto&y_{1-i}\\
y_i&\longmapsto&x_{1-i}
\end{array}.
\end{equation}
By considering the map $\varphi\circ h\circ\varphi^+$,
we deduce the existence of a graded algebra morphism
\begin{equation}\label{lede}
\begin{array}{rcl}
h_+~:\quad \CC[S_{\infty}\backslash B_-]_q^+&\longrightarrow&A_q\\
m(\lambda)&\longmapsto&\lim\limits_{N\to\infty}
\dfrac{x_1^{-1}}{1+\dfrac{(\lambda x_1y_1)^{-1}}
       {1+\dfrac{(\lambda y_{1}x_2)^{-1}}{\dfrac{\ddots
       (\lambda y_{N-1}x_{N})^{-1}}
       {1+(\lambda x_{N}y_{N})^{-1}}}}}
\end{array}
\end{equation}
with the convention that $\dfrac{a}{b}=ab^{-1}$.
Explicitly, thanks to Proposition \ref{aba},
we get
\begin{equation}\label{defmiq}
     h_+(m_i)=\sum F_q(\a_1,\a_2,\ldots)\,
     x_1^{-1}(x_1 y_1)^{-\a_1}\ldots
     (x_{k}y_{k})^{-\a_{2k-1}}
     (y_{k}x_{k+1})^{-\a_{2k}}\ldots
\end{equation} 
As usual, the sum is taken on all almost zero sequences $(a_k)$ 
such that $\sum_k a_k=i-1$.
Here also, in each new term of the sequence 
$m_i$ occurs one and only one new variable of the form
$x_k$ or $y_k$ (according to the parity of $i$).
Therefore, the same argument as before shows the two following results.
\begin{prop}\label{unpo}
A basis for $\CC[S_{\infty}\backslash B_-]_q^+$ 
is given by the family 
$\eta_a:=\prod_{i=1}^{\infty}m_i^{a_i}$ where 
$a=(a_i)_{i\in\NN^*}$ denotes any almost zero sequence of integers.
\end{prop}

\begin{cor}
The map $h_+$ est injective.
\end{cor}
 
In the classical case, $h_+$ is also a birational isomorphism 
from $\CC[S_{\infty}\backslash B_-]^+$ to
$\CC[x_i^{-1},y_i^{-1},\, i>0]$. 

Let us consider again the quantum algebra 
$\CC[H_-\backslash B_-]_q$ defined in Section
\ref{resprin}.

\begin{lem}\label{doub}
The natural map
$$
\begin{array}{rcl}
\CC[H_-\backslash B_-]_q&\longrightarrow&\CC[S_{\infty}\backslash B_-]_q
\ootimes
\CC[S_{\infty}\backslash B_-]_q^+\\
\forall\, i,\quad u_i&\longmapsto&u_i\ootimes 1\\
m_i&\longmapsto&1\ootimes m_i
\end{array}
$$
is a graded algebra isomorphism.
\end{lem}

\begin{dem}
It is enough to construct the inverse.
But, there exist natural morphisms $f$ and $f^+$ from 
$\CC[S_{\infty}\backslash B_-]_q$ and $\CC[S_{\infty}\backslash B_-]_q^+$
to $\CC[H_-\backslash B_-]_q$. These morphisms $q^{-1}$-commute.
So, there exists $f\ootimes f^+:\;
\CC[S_{\infty}\backslash B_-]_q
\ootimes\CC[S_{\infty}\backslash B_-]_q^+\longrightarrow
\CC[H_-\backslash B_-]_q$. One can check that 
this gives an inverse for the studied map.
\end{dem}

Hence, by virtue of Proposition \ref{deff} and \ref{unpo}, 
we deduce that $\CC[H_-\backslash B_-]_q$
is indeed a flat deformation of $\CC[H_-\backslash B_-]$.
This completes the proof of Proposition \ref{thb}.
Note that $\CC[H_-\backslash B_-]_q$ is a graded algebra
with the gradation given by $\deg u_i=-\deg m_i$
for all $i$ and that a basis of $\CC[H_-\backslash B_-]_q$
is given by the family
$\prod_{i=1}^{\infty} u_i^{\alpha_i}\prod_{j=1}^{\infty} 
m_j^{\beta_j}$ where $(\alpha_i)$
and $(\beta_j)$ are two almost zero sequences of integers.

\section{Quantum Drinfeld-Sokolov correspondence}
The aim of this section is to prove Theorem \ref{thc}.
As a result of the previous section, 
we have already proved the existence of $\DS_q$.
Indeed, it suffices to consider the morphisms $h$ and $h_+$
seen in (\ref{morh}) et (\ref{lede}), to note that
$\Im h\subset\CC[x_i^{\pm 1},y_i^{\pm 1},\, i\leq 0]_q$,
$\Im h_+\subset\CC[x_i^{\pm 1},y_i^{\pm 1},\, i>0]_q$
and to take into account Lemma \ref{doub} together with
the isomorphism $A_q\simeq \CC[x_i^{\pm 1},y_i^{\pm 1},\, i\leq 0]_q\ootimes
\CC[x_i^{\pm 1},y_i^{\pm 1},\, i>0]_q$.
Note that the injectivity of $h$ and of $h_+$ imply
the one of $\DS_q$. It remains to prove the equality
$\DS_q\circ H_{\mu}=\ad(I_{\mu})\circ\DS_q$.
For that, the idea is first to prove the existence of
$H_n$ (this will be achieved in subsection \ref{exx})
and then, using the embedding 
of $\CC[H_-\backslash B_-]_q$ into $A_q$,
to extend $H_n$ not only on $A_q$ but
also on ${\bar A}_q$ the algebra obtained from $A_q$
by adding the two half screening charges $\SPMB$ of 
the discrete sine-Gordon system. The interest in considering
this algebra is that it is endowed with an $U_q\bm$-module-algebra structure,
where $\bm$ is a Borel subalgebra of $\widehat{{\mathfrak{sl}}_2}$.
Moreover, the adjoint actions of integrals of motion extend
to ${\bar A}_q$ and commute with the action of $U_q\bm$.
Conversely, each derivation which commutes with the action
of $U_q\bm$ is the adjoint action of an integral of motion.
We shall use this fact to complete the proof.
First, we start by giving precise definitions of the quantum group
$U_q\bm$ and algebras ${\bar A}_q$ 
and $\CC[H_-\backslash B_-N_+]_q$.

\begin{Def}
{\it{Let $\bm$ be a Borel subalgebra of $\widehat{{\mathfrak{sl}}_2}$.
We note $U_q\bm$ the quantum group given by generators~:
$k_{\varepsilon}^{\pm 1}$, 
$f_{\varepsilon'}$, 
with $\varepsilon,\varepsilon'\in\{ +,-\}$
and relations~:
\begin{align}
k_{\varepsilon}k_{\varepsilon}^{-1}&=
k_{\varepsilon}^{-1}k_{\varepsilon}=1\\
k_{\varepsilon}k_{\varepsilon'}&=k_{\varepsilon'}k_{\varepsilon}\\
k_{\varepsilon}f_{\varepsilon'}k_{\varepsilon}^{-1}&=
q^{\a_{\varepsilon,\varepsilon'}}f_{\varepsilon'},
\end{align}
together with the quantum Serre relations between
$f_{\pm}$ and $f_{\mp}$:
\begin{equation}
f_{\pm}^{3}f_{\mp} -(q+1+q^{-1})\bigl( f_{\pm}^{2}f_{\mp}f_{\pm}-
f_{\pm}f_{\mp}f_{\pm}^{2}\bigr)-f_{\mp}f_{\pm}^{3}=0,
\end{equation}
with the convention that $\a_{\varepsilon,\varepsilon'}=1$ if 
$\varepsilon=\varepsilon',\, -1$ if not.
The comultiplication is given by~:
\begin{align}
   \Delta(f_{\varepsilon})&=f_{\varepsilon}\otimes 1
   +k_{\varepsilon}\otimes f_{\varepsilon}\\
   \Delta(k_{\varepsilon}^{\pm 1})&=k_{\varepsilon}^{\pm 1}\otimes
    k_{\varepsilon}^{\pm 1}.
\end{align}
We shall use neither the antipode nor the co-unit in this article.}}
\end{Def}
 
\subsection{The extended phase space ${{\bar{A}}}_q$.}\label{eten}
If we consider only a finite number of sites 
$x_i^{\pm 1},y_i^{\pm 1},\, i\in\{ 1,\ldots,n\}$,
it can be shown from Lemma \ref{enrf} that we get a $U_q\bm$-module-algebra.
For all $x$ homogeneous with respect to $\deg$, 
the formulas are the following~:
\begin{align}
f_{+}.x&=\sum_{i=1}^n[x_i,x]_q,\\
f_{-}.x&=\sum_{i=1}^n[y_i,x]_q,\\
k_{\pm}.x&=q^{\pm\deg{x}}x.
\end{align}
If we consider now an infinite number of sites at the left 
of an arbitrary site $x_i^{\pm 1},y_i^{\pm 1},\, i\leq N$, 
we also obtain a $U_q\bm$-module-algebra. For that, we set~:
\begin{align}
\label{poursp}
f_{+}.x&=\sum_{i\leq N}[x_i,x]_q,\\
\label{poursm}
f_{-}.x&=\sum_{i\leq N}[y_i,x]_q,\\
\label{pourk}
k_{\pm}.x&=q^{\pm\deg{x}}x,
\end{align}
for any $x$ homogeneous with respect to $\deg$.
This follows from the fact that for all $x\in A_q$,
$x_i$ and $x$ $q$-commute
provided that $i$ is small enough.
However, if we consider the whole algebra $A_q$, there is no longer
a $U_q\bm$-module-algebra structure on it. For that, it is 
necessary to add the half screening charges $\SPB$ and $\SMB$
which correspond heuristically to 
$\sum_{i>0}x_i$ and $\sum_{i>0}y_i$.

\begin{Def}
{\it{We note ${{\bar{A}}}_q$
the algebra given by generators~:
$\SPB$. $\SMB$, $x_i^{\pm 1}$, $y_i^{\pm 1}$, $i\in\ZZ$ 
and relations~: 
\begin{equation}\label{defaqb}
   \begin{array}{rrcl}
     \forall i<j,\quad &x_ix_j&=&qx_jx_i\\
     &y_iy_j&=&qy_jy_i\\
     &y_ix_j&=&q^{-1}x_jy_i\\
     \forall i\leq j,\quad &x_iy_j&=&q^{-1}y_jx_i\\
     \forall i\in\ZZ,\quad &x_i\SPB-q\SPB x_i
     &=&\sum_{j=1}^{i}[x_i,x_j]_q\\
     &x_i\SMB-q^{-1}\SMB
     &=&\sum_{j=1}^{i}[x_i,y_j]_q\\
     &y_i\SPB-q^{-1}\SPB y_i
     &=&\sum_{j=1}^{i}[y_i,x_j]_q\\
     &y_i\SMB-q\SMB y_i
     &=&\sum_{j=1}^{i}[y_i,y_j]_q
   \end{array}
\end{equation}
together with the quantum Serre relations
between $\SPMB$ and $\SMPB$~:
  \begin{equation}\label{serreqspm}
   \SPMB^{3}\SMPB -(q+1+q^{-1})\bigl( \SPMB^{2}\SMPB\SPMB-
   \SPMB\SMPB\SPMB^{2}\bigr)-\SMPB\SPMB^{3}=0.
  \end{equation}
As usual, for two elements $a$ and $b$,
$[a,b]_q$ denotes the $q$-commutator of $a$ and $b$. 
The gradation $\deg$ on 
${{\bar{A}}}_q$ is given by~:
$$
\forall\, i\in\ZZ,\quad \deg{x_i}=\deg{\SPB}=-\deg{y_i}=-\deg{\SMB}=1.
$$}}
\end{Def}

The following result can be proved easily.
\begin{lem}\label{kouch}
A basis for ${\bar A}_q$ is given by the family
$\prod_{i=-\infty}^{+\infty}x_i^{\a_i}y_i^{\b_i}u$
where $u$ belongs to a basis ${\cal B}$ of
$\CC[\SPB,\SMB]_q\simeq U_q\nm$ and
$(\a_i)$, $(\b_i)$ are two almost zero sequences in
${\ZZ}^{\ZZ}$.
\end{lem}
Hence, thanks to subsection \ref{defprin},
we get the following lemma.

\begin{lem}
There is a natural graded algebra embedding 
$A_q\hookrightarrow {\bar A}_q$. 
This embedding identifies generators 
$x_i^{\pm 1}$ and $y_i^{\pm 1},\, i\in\ZZ$ 
with the ones of ${\bar A}_q$.
\end{lem}

The semi-classical limit of ${\bar A}_q$ is 
$\CC\bigl[x_i^{\pm 1},y_i^{\pm 1}, i\in\ZZ,\,
\SEPAB,
\{\SEPAB,\SEPBB\},
\{\SEPAB,\{\SEPAB,\SEPBB\}\},\ldots,\varepsilon_k\in\{+,-\}\bigr]$.
We note this algebra ${{\bar{A}}}_{\cl}$.
Let us remark that it is possible to extend the half-translation
automorphism $T^{\frac{1}{2}}$ on ${{\bar{A}}}_{\cl}$
by setting $T^{\frac{1}{2}}\bigl(\SPB\bigr)=\SMB$
and $T^{\frac{1}{2}}\bigl(\SMB\bigr)=\SPB-x_1$.
It can be shown that ${{\bar{A}}}_{\cl}$ is the localized 
of a subalgebra of a projective limit of algebras.
Explicitly, these algebras are the ones generated by
variables $x_i$ and $y_i$ for $i\leq n$ with obvious
projection morphisms. The considered subalgebra is the one generated
by $x_i$, $y_i$, $i\in\ZZ$ and the half-screening charges $\SPB$
and $\SMB$ identified with $(x_1,x_1+x_2,\ldots)$
and $(y_1,y_1+y_2,\ldots)$. The multiplicative set
is generated by elements $x_i^{-1}$ and $y_i^{-1}$
for $i\in\ZZ$. It satisfies the Ore relation
\cite{Dix}. Therefore, we can deduce from formulas 
(\ref{poursp}),
(\ref{poursm}), (\ref{pourk}) 
that there exists a $U_q\bm$-module-algebra structure 
on ${{\bar{A}}}_q$ given by~:
\begin{align}
\label{fepsx1}
f_{+}.x&=\bigl[\SP,x\bigr]_q=
\sum_{i=-\infty}^{0}[x_i,x]_q
+\bigl[\SPB,x\bigr]_q\\
\label{fepsx2}
f_{-}.x&=\bigl[\SM,x\bigr]_q:=
\sum_{i=-\infty}^{0}[y_i,x]_q
+\bigl[\SMB,x\bigr]_q\\
\label{kpmx}
k_{\pm}.x&=q^{\pm\deg{x}}x.
\end{align}
At the semi-classical limit, it also gives 
a $U\bm$-module-algebra structure 
on ${{\bar{A}}}_{\cl}$.
\subsection{The quantum homogeneous space
$\CC[H_-\backslash B_-N_+]_q$}\label{avecs}
Geometrically, at the classical level, adding 
half screening charges is the same as
studying the Schubert cell $B_-N_+$ of $G$ instead of its
Borel group $B_-$.

\begin{Def}
{\it{We denote by $\CC[H_-\backslash B_-N_+]_q$ 
the quantum algebra given by generators~:
$\SPB$, $\SMB$, $u_i,\, m_i,\, i\in\ZZ$ and relations~:
\begin{gather}
\label{peru}
\bigl(\lambda^{-1}u(\lambda)-\mu^{-1}u(\mu)\bigr)\bigl(u(\lambda)-u(\mu)\bigr)
=q\bigl(u(\lambda)-u(\mu)\bigr)
\bigl(\lambda^{-1}u(\lambda)-\mu^{-1}u(\mu)\bigr)\\
\label{perm}
\bigl(\lambda^{-1}m(\lambda)-\mu^{-1}m(\mu)\bigr)\bigl(m(\lambda)-m(\mu)\bigr)
=q^{-1}\bigl(m(\lambda)-m(\mu)\bigr)
\bigl(\lambda^{-1}m(\lambda)-\mu^{-1}m(\mu)\bigr)\\
\label{perum}
u(\lambda)m(\mu)=q^{-1}m(\mu)u(\lambda)\\
\label{perus}
u(\lambda)\SPMB=q^{\pm 1}\SPMB u(\lambda)\\
\label{permsp}
m(\lambda)\SPB-q^{-1}\SPB m(\lambda)=1-q^{-1}\\
\label{permsm}
m(\lambda)\SMB-q\SMB m(\lambda)=-(q-1)\lambda^{-1}m(\lambda)^{2}\\
\label{but}
\SPMB^3 \SMPB-(q+1+q^{-1})\bigl(\SPMB^2\SMPB\SPMB
-\SPMB\SMPB\SPMB^2\bigr)-\SMPB\SPMB^3=0
\end{gather}
with the same notation as before i.e.,
$u(\lambda)=\sum_{i=0}^{\infty} (-1)^i u_{i+1}\lambda^{-i}$
and
$m(\lambda)=\sum_{i=0}^{\infty} (-1)^i m_{i+1}\lambda^{-i}$.}}
\end{Def}

The algebra $\CC[H_-\backslash B_-N_+]_q$ 
is graded with the gradation given by~:
$$
\forall\, i\in\ZZ,\quad \deg{u_i}=\deg{\SPB}=-\deg{m_i}=-\deg{\SMB}=1.
$$
The relations between $u(\lambda)$, $m(\lambda)$ and $\SPMB$ will 
appear to be natural when we prove the following proposition
which claims the existence of the morphism $\DSB_q$.
\begin{prop}\label{prolonm}
The map 
$$
\begin{array}{rcl}
\DSB_q~:\quad\CC[H_-\backslash B_-N_+]_q&\longrightarrow&{{\bar{A}}}_q\\
\forall\, i\in\NN^*,\quad u_i&\longmapsto&\DS_q(u_i)\\
m_i&\longmapsto&\DS_q(m_i)\\
\SPMB&\longmapsto&\SPMB
\end{array}
$$
exists and defines a graded algebra morphism.
\end{prop}

\begin{dem}
We need to prove some compatibility relations.
The ones dealing with $\DS_q(u_k)$, $k\in\NN^*$ and
$\SPMB$ follow from the fact that all the terms in
$\DS_q(u_k)$ are sums and products of
$x_i^{\pm 1}$ and $y_i^{\pm 1}$ for $i\leq 0$.
The ones dealing with $\DS_q(m_k)$ and $\SPB$
can be handled in the following way.
If we take again the involution $\varphi$ defined in
(\ref{definfi}), we have $\DS_q(m_k)=\varphi\bigl(\DS_q(u_k)\bigr)$. 
So, for any integer $n$ large enough,
\begin{align*}
\bigl[\DS_q(m_k),\SPB\bigr]_q&=\bigl[\DS_q(m_k),\sum_{i=1}^n x_i\bigr]_q=
\varphi\Bigl(\bigl[\sum_{i=-n+1}^0 y_i,\DS_q(u_k)\bigr]_q\Bigr)\\
&=(T^{n}\circ\varphi)\Bigl
(\bigl[\sum_{i=1}^{n} y_i,T^{n}\bigl(\DS_q(u_k)\bigr)\bigr]_q\Bigr).
\end{align*}
Recall the notation of subsection \ref{cell}
and in particular the morphism
$h_n:\CC[B_- w_n B_-\cap N_+]_q\longrightarrow A_q$,
we have
$T^{n}\DS_q(u_k)=h_n(\alpha_k)$ and
$h_n(a_{1,2}^{(0)})=\sum_{i=1}^{n}y_i$.
Then the result comes from commuting relations in
$\CC[B_- w_n B_-\cap N_+]_q$.
We prove the compatibility relation between $\DS_q(m_k)$
and $\SMB$ using a similar method. 
\end{dem}

{}From the commuting relations in $\CC[H_-\backslash B_-N_+]_q$
together with the results of section \ref{secquanti},
we can deduce the following corollary.

\begin{cor}\label{defobn}
If ${\cal B}$ denotes a basis for $\CC[\SMB,\SPB]_q\simeq U_q\nm$,
then the family $\prod_{i=1}^{\infty}u_i^{\a_i}\prod_{j=1}^{\infty}
m_j^{\beta_j}u$
where $(\a_i)$ and $(\beta_j)$
are two almost zero sequences of integers and 
$u\in{\cal B}$ is a basis for $\CC[H_-\backslash B_-N_+]_q$.
\end{cor}

We also obtain the following result.
\begin{cor}
The morphism $\DSB_q$ is injective. 
\end{cor}

\begin{dem}
This is a consequence of Corollary \ref{defobn},
Corollary \ref{kouch}
and of the already seen fact that
each new term of the sequence $u_k$ (resp. $m_k$)
gives a new variable $x_{-i}$ (resp. $x_{i}$)
or $y_{-i}$ (resp. $y_{i}$), $i\in\NN$ according to the 
parity of $k$.
\end{dem}

Corollary  \ref{defobn} also shows that 
$\CC[H_-\backslash B_-N_+]_q$ is a flat deformation 
of the function algebra of the Poisson manifold
$\bigl(H_-\backslash B_-N_+,P_{\infty}\bigr)$.
Poisson relations on this manifold show that we obtain
a quantum deformation. Moreover, the map~:
\begin{equation}\label{sect}
   \begin{array}{rcl}
F~:\quad\CC[[\lambda^{-1}]]\times\lambda^{-1}\CC[[\lambda^{-1}]]
\times N_+&\longrightarrow&
    \bigl(H_-\backslash B_-,N_+\bigr)\\
     (u_{\cl}(\lambda),v_{\cl}(\lambda),n_+)&\longmapsto&\Bigl(\cl_{H_-} 
     \begin{pmatrix}
      1&v_{\cl}(\lambda)\\
      0&1 
     \end{pmatrix}
     \begin{pmatrix}
      1&0\\
     u_{\cl}(\lambda)&1
     \end{pmatrix},n_+\Bigr)
    \end{array}
  \end{equation}
is a bijection and the elements $\SPB$ and $\SMB$
correspond classically to the functions 
$a_{2,1}^{(1)}$ and $a_{1,2}^{(0)}$ on $N_+$.
On the other hand, by virtue of subsection \ref{ahbon},
classical limits of $u_i$ and $m_i$ correspond to coordinate
functions with generating functions 
$u_{\cl}(\lambda)$ and $m_{\cl}(\lambda)$.
Note that Corollary \ref{defobn}
also implies the following lemma.
\begin{lem}\label{plon}
There is a natural graded algebra embedding
$\CC[H_-\backslash B_-]_q\hookrightarrow\CC[H_-\backslash B_-N_+]_q$.
This embedding identifies generators $u_i$ (resp. $m_i$), $i\in\ZZ$
of $\CC[H_-\backslash B_-]_q$ with the ones of
$\CC[H_-\backslash B_-N_+]_q$. 
Moreover, we have the following commutative diagram 
where all maps are graded algebras embedding~:
$$
\begin{array}{ccc}
 \CC[H_-\backslash B_-]_q&\hookrightarrow&\CC[H_-\backslash B_-N_+]_q\\
 \DS_q\downarrow&&\downarrow\DSB_q\\
 A_q&\hookrightarrow&{\bar{A}}_q.
\end{array}
$$
\end{lem}
Proposition \ref{prolonm}
together with Lemma \ref{plon}
lead to the following result.

\begin{cor}\label{stru} 
There is a $U_q\bm$-module-algebra structure on 
$\CC[H_-\backslash B_-N_+]_q$ given by~:
\begin{align}
     \label{zaz}
     f_+.u(\lambda)&=-(q-1)\lambda^{-1}u(\lambda)^{2}-(q-q^{-1})
     u(\lambda)\SPB\\
     f_-.u(\lambda)&=(1-q^{-1})+(q-q^{-1})u(\lambda)\SMB\\
     k_{\pm}.u(\lambda)&=q^{\pm 1}u(\lambda)\\
     f_+.m(\lambda)&=1-q+(q-q^{-1})m(\lambda)\SPB\\
     f_{-}.m(\lambda)&=(1-q^{-1})\lambda^{-1}m(\lambda)^{2}-(q-q^{-1})
     m(\lambda)\SMB\\
     k_{\pm}.m(\lambda)&=q^{\mp 1}m(\lambda)\\
     f_{\varepsilon}.\SEPPB&=\bigl[\SEPB,\SEPPB\bigr]_q.
\end{align}
The morphism $\DSB_q$ defined above is a 
$U_q\bm$-module-algebra morphism.
\end{cor}

\begin{dem}
This comes from the fact that the morphism $\DSB_q$ is injective
and from the computation of 
$f_{\pm}.\DSB_q(x)$ for $x\in\CC[H_-\backslash B_-N_+]_q$.
For example, according to formula (\ref{fepsx1}),
the computation of $f^+.\DSB_q(u_k)$ leads to the computation of 
$\sum_{i=-\infty}^0 \bigl[x_i,\DSB_q(u_k)\bigr]_q$. 
This sum is finite. So, recalling the involution $\varphi$ on $A_q$
defined in (\ref{definfi}), we have~:
\begin{align*}
\sum_{i=-\infty}^0 \bigl[x_i,\DSB_q(u_k)\bigr]_q&=
\sum_{i=1}^{\infty}\varphi\Bigl(\bigl[\DSB_q(m_k),y_i\bigr]_q\Bigr)=
\varphi\Bigl(\bigl[\DSB_q(m_k),\SMB\bigr]_q\Bigr)\\
&=(\varphi\circ\DSB_q)\Bigl(\bigl[m_k,\SMB\bigr]_q\Bigr).
\end{align*} 
Then it suffices to use (\ref{permsm}) to get the expression of
$f^+.\DSB_q(u_k)$.
\end{dem}

\subsection{Adjoint action of integrals of motion on 
${{\bar{A}}}_q$}\label{acad}
Let $I$ be an integral of motion. By using the definition 
seen in (\ref{adjoin}) of the adjoint action of $I$ on $A_q$
together with the equality between (\ref{alij1}) and (\ref{alij2}), 
it can be shown that there exists an unique homogeneous element 
$R_+(I)\in A_q[1]$ without constant term such that 
$\ad(I)(x_1)=T\bigl(R_+(I)\bigr)-R_+(I)$.
For instance, if $I=I_1$ is the first integral of motion 
with respect to the basis $(I_k)$ of ${\cal I}$ seen in
Proposition \ref{basi}, then $R_+(I_1)=-y_0^{-1}$.
It follows that for any integer $n$, we have
$\ad(I)(x_1+\ldots+x_n)=T^n\bigl(R_+(I)\bigr)-R_+(I)$.
On the other hand, thanks to the form taken by the $I_k$, 
we have $\ad(I)\circ T^{\frac{1}{2}}=T^{\frac{1}{2}}\circ\ad(I)$
for all $I\in{\cal I}$. So, there exists also
$R_-(I)\in A_q[-1]$ with
$R_-(I)=T^{\frac{1}{2}}\bigl(R_+(I)\bigr)$ such that
$\ad(I)(y_1)=T\bigl(R_-(I)\bigr)-R_-(I)$.
So, for all $n$, $\ad(I)(y_1+\ldots+y_n)=T^n\bigl(R_-(I)\bigr)-R_-(I)$.
This leads to extend the derivation $\ad(I)$ on 
${{\bar{A}}}_q$ as explained in the following proposition.
\begin{prop}\label{onv}
Let $\Der({{\bar{A}}}_q)$ be the Lie algebra of  
derivations on ${\bar A}_q$.
For $I\in{\cal I}$, there is a unique derivation $\ad(I)$
on ${{\bar{A}}}_q$ which satisfies
formula (\ref{adjoin}) if $x$ belongs to
$A_q$ and $\ad(I)(\SPMB)=-R_{\pm}(I)$.
Moreover, the kernel of the Lie algebra morphism~:
\begin{equation}\label{micoud}
\begin{array}{rcl}
\ad~:\quad{\cal I}&\longrightarrow&\Der({{\bar{A}}}_q)\\
I&\longmapsto&\ad(I)
\end{array}
\end{equation}
is $\CC[q,q^{-1}]$ i.e., the one-dimensional Lie subalgebra of all
constant integrals of motion. Its image is 
$\Der_{U_q\bm}({{\bar{A}}}_q)$
the Lie subalgebra of all derivations which commute
with the action of $U_q\bm$.
\end{prop}

\begin{dem}
Let $I\in{\cal I}$. So as to prove the existence of
$\ad(I)$ on ${{\bar{A}}}_q$, it is necessary to
show some compatibility relations like~:
\begin{equation}\label{compati}
   \forall\, j\in\ZZ,\quad
   \bigl[ \ad(I)(x_j),\SPB\bigr]_q
   +\bigl[ x_j,-R_{+}\bigr]_q
   =\sum_{k=1}^{j}\ad(I)\bigl( [x_j,x_k]_q\bigr).
 \end{equation}
But, according to Proposition \ref{nondemb},
$\ad(I)$ is well defined on $A_q$.
So, for any fixed integer $j$, we have~:
$$
\forall\, n\in\NN,\quad
   \bigl[ \ad(I)(x_j),\SPN\bigr]_q
   +\bigl[ x_j,T^n (R_{+})-R_{+}\bigr]_q
   =\sum_{k=1}^{j}\ad(I)\bigl( [x_j,x_k]_q\bigr)
$$
with $\SPN =\sum_{k=1}^n x_k$. 
This equality leads to (\ref{compati})
by taking $n$ large enough.
The other relations except those coming from the quantum Serre relations
between $\SPMB$ and $\SMPB$ can be proved in the same way.
The unicity of $\ad(I)$ is obvious.

To prove that $\ad(I)$ and $f_{\pm}$ commute, we set
$C_{\pm}=\{ x\in{{\bar{A}}}_q/\,
\ad(I)\circ f_{\pm}(x)=f_{\pm}\circ\ad(I)(x)\}$.
We remark that for all $x$ homogeneous with respect to
$\deg$, we have $\deg\bigl(\ad(I)(x)\bigr)=\deg(x)$.
Hence, $C_{\pm}$ is a graded subalgebra of ${{\bar{A}}}_q$.
Then, computation shows that 
$A_q\cup\{ \SPB,\SMB\}\in C_{\pm}$. It follows that
$C_{\pm}={{\bar{A}}}_q$.

Conversely, let us fix $D\in\Der_{U_q\bm}({{\bar{A}}}_q)$. 
Using the fact that the result is 
true at the classical level (\cite{ENR}) and the fact that classical
integrals of motions can be quantized, first, we show the following result~:
\begin{equation}\label{ais}
\forall \delta\in\Der_{U_q\bm}({{\bar{A}}}_q)\, \forall n\in\NN\,
\exists\, I\in{\cal I}\, \exists\, \delta'\in\Der_{U_q\bm}({{\bar{A}}}_q),\,
\forall x\in{{\bar{A}}}_q,\;
\delta(x)=\ad(I)(x)+(q-1)^n \delta'(x).
\end{equation}
Then, we can deduce from (\ref{ais}) and from
the explicit form of the basis $(I_k)$
of ${\cal I}$ that (1) 
$D$ is a graded derivation (it means that if $x\in{{\bar{A}}}_q$
is homogeneous with respect to $\deg$, then $D(x)$
is also homogeneous with respect to $\deg$ and $\deg D(x)=\deg x$), (2) 
$[D,T^{\frac{1}{2}}]=0$ or in other words,
$D$ and $T^{\frac{1}{2}}$ commute,
and (3) $A_q$ is invariant by $D$.
Hence, we show that $D$ is entirely defined by $D(x_0)$~:
the natural map coming from the foregoing,
\begin{equation}\label{injder}
\begin{array}{rcl}
\Der_{U_q\bm}({{\bar{A}}}_q)&\longrightarrow&\Der_{T^{\frac{1}{2}}}(A_q)\\
\delta&\longmapsto&\delta_{|A_q}
\end{array}
\end{equation}
is injective. 
Let us give some definitions.
Let $V_p$ be the free sub-module of $A_q$ of all homogeneous elements
of degree $p$ with respect to the principal gradation $\deg_p$
(see subsection \ref{defgra}).
Let ${\cal B}_p$ be a basis for $V_p$. If there exists $q\in\ZZ$
such that $p=-2q+1$, then ${\cal B}_p$ can be chosen such that 
$\ad(I_q)(x_0)$ belongs to ${\cal B}_p$. There is $N\in\NN$ and
$\a_p\in V_p$, $p\in\{-N,\ldots,N\}$ such that 
$D(x_0)=\sum_{p=-N}^{N}\a_p$. Let $p\in\{-N,\ldots,N\}$.
By projecting (\ref{ais}) on $V_p$ with $\delta=D$ and $x=x_0$, 
we see that the valuations in $q-1$ of all coefficients of 
$\a_p$ on basis ${\cal B}_p$
(except perhaps the element $\ad(I_q)(x_0)$ of ${\cal B}_p$
if there is an integer $q$ such that $p=-2q+1$)
are arbitrarily large. Hence, $\a_p=0$ or $\a_p$ is proportional
to $\ad(I_q)(x_0)$. 
Thus, there is $I\in{\cal I}$ such that 
$D(x_0)=\ad(I)(x_0)$. Then, the injectivity of the map
(\ref{injder}) shows that $D=\ad(I)$.
\end{dem}
\subsection{Existence of $H_n$}\label{exx}
The aim of this section is to prove Proposition
\ref{exih}.
In fact, we shall prove the existence of $H_n$ not only
on $\CC[H_-\backslash B_-]_q$ but also on 
$\CC[H_-\backslash B_-N_+]_q$. This will imply Proposition 
\ref{exih}.
\begin{prop}\label{proex}
There is a commutative family of derivations 
$(H_n)_{n\in\NN^*}$
on $\CC[H_-\backslash B_-N_+]_q$ 
which quantizes the classical action by vector fields of $\fh_+$ on 
$H_-\backslash B_-N_+$ and which commute with the action of
$U_q\bm$ on $\CC[H_-\backslash B_-N_+]_q$
(for the definition of $\fh_+$, see subsection \ref{defprin}).
If $H(\mu)=\sum_{k=1}^{\infty}(-1)^k H_k \mu^{-k}$
denotes the generating function of $(H_n)_{n\in\NN^*}$,
then, the derivations $H_n$ are defined by formulas~:
\begin{align}
\label{jmul}
H_{\mu}\bigl(u(\lambda)\bigr)&=\dfrac{1}{\lambda^{-1}-\mu^{-1}}
\bigl(\lambda^{-1}u(\lambda)-\mu^{-1}u(\mu)\bigr)v(\mu)u(\lambda)
-\dfrac{\mu^{-1}}{\lambda^{-1}-\mu^{-1}}
\bigl(u(\lambda)-u(\mu)\bigr)\bigl(1+v(\mu)u(\mu)\bigr)\\
\label{jmum}
H_{\mu}\bigl(m(\lambda)\bigr)&
=\dfrac{\mu^{-1}}{\lambda^{-1}-\mu^{-1}}\bigl(1+m(\mu)w(\mu)\bigr)
\bigl(m(\lambda)-m(\mu)\bigr)-\dfrac{1}{\lambda^{-1}-\mu^{-1}}m(\lambda)w(\mu)
\bigl(\lambda^{-1}m(\lambda)-\mu^{-1}m(\mu)\bigr)\\
H(\mu)(\SPB)&=-\mu^{-1}\bigl(u(\mu)+u(\mu)v(\mu)u(\mu)\bigr)\\
H(\mu)(\SMB)&=v(\mu)
\end{align}
with 
$v(\mu)=-{\bigl(u(\mu)+\mu m(\mu)^{-1}\bigr)}^{-1}$ and
$w(\mu)=-{\bigl(m(\mu)+\mu u(\mu)^{-1}\bigr)}^{-1}$.
\end{prop}
\subsection{The classical case}
For $n\in\NN$, set $h_n= \diag(\lambda^n,-\lambda^n)$.
With the notation of subsection \ref{avecs}
and in particular of the map $F$ defined in (\ref{sect}),
we show that the left translation action of $h_n$
on generating series $u_{\cl}(\lambda)$ and
$v_{\cl}(\lambda)$ of coordinate functions $u_i$ and
$m_i$ is given by formulas~:
\begin{align}
\label{huc}
h_n.u_{\cl}&=\bigl[ (1+2u_{\cl}v_{\cl})\lambda^{n}\bigr]_{\leq}u_{\cl}
-2\bigl[ u_{\cl}(1+u_{\cl}v_{\cl})\lambda^{n}\bigr]_{\leq}- \bigl[
2v_{\cl}\lambda^{n}\bigr]_{<}u_{\cl}^{2}+\bigl[(1+2u_{\cl}v_{\cl})
\lambda^{n}\bigr]_{\leq}u_{\cl}\\
\label{hvc}
h_n.v_{\cl}&=2\bigl[v_{\cl}\lambda^{n}\bigr]_{<}(1+2u_{\cl}v_{\cl})
-2\bigl[(1+2u_{\cl}v_{\cl})\lambda^{n}\bigr]_{\leq} 
\end{align}
where for $x\in\CC((\lambda^{-1}))$,
$x_{\leq}$ (resp. $x_{<}$) denotes the part of
$x$ in $\CC[[\lambda^{-1}]]$ (resp. $\lambda^{-1}\CC[[\lambda^{-1}]]$).
Let $h(\mu)$ be the generating series of $h_n$. 
{}From (\ref{huc}) and (\ref{hvc}), 
it is possible to compute the action of 
$\dfrac{1}{2}h(\mu)$ on $u_{\cl}(\lambda)$ and $m_{\cl}(\lambda)$.
It can be checked that these relations are precisely the ones we get from
(\ref{jmul}) and (\ref{jmum}) when $q\to 1$.
In the same way, with the identification of $\SPB$
and $\SMB$ with $1\otimes a_{2,1}^{(1)}$ 
and $1\otimes a_{1,2}^{(0)}$, it can be shown that we have
$\dfrac{1}{2}h(\mu).\SPB=-\mu^{-1}\bigl( u(\mu)+u(\mu)^2v(\mu)\bigr)$
and $\dfrac{1}{2}h(\mu).\SMB=v(\mu)$. 
Hence, it is clear that if derivation $H_n$ exists then it deforms
the classical action of $\dfrac{1}{2}h_n$ on the homogeneous space
$H_-\backslash B_-N_+$.

\subsection{The algebra $U_{\lambda,\mu}$}
To obtain algebra $U_{\lambda,\mu}$, we just need to replace generating series
$u(\lambda)$ and $u(\mu)$ by variables $u_{\lambda}$
and $u_{\mu}$.

\begin{Def}
{\it{We denote by $U_{\lambda,\mu}$ the algebra over the ground ring
$\CC[\lambda^{-1},\mu^{-1}]$ given by generators
$u_{\lambda}$ and $u_{\mu}$ and relation~:
$\bigl(\lambda^{-1}u_{\lambda}-\mu^{-1}u_{\mu}\bigr)
\bigl(u_{\lambda}-u_{\mu}\bigr)
=q\bigl(u_{\lambda}-u_{\mu}\bigr)
\bigl(\lambda^{-1}u_{\lambda}-\mu^{-1}u_{\mu}\bigr)$.}}
\end{Def}

Note that $U_{\lambda,\mu}$ is a graded algebra with 
respect to the gradation given by $\deg u_{\lambda}=
\deg u_{\mu}=1$. It can be proved that $U_{\lambda,\mu}$
does not have any torsion of zero divisor. and that a basis
is given by the family $(u_{\lambda}^{\a}u_{\mu}^{\b})$ 
with $(\a,\b)\in\NN^2$.
Thanks to definition of $U_{\lambda,\mu}$ and 
$\CC[S_{\infty}\backslash B_-]_q$, there is an algebra morphism~:
\begin{equation}\label{morspeu}
\begin{array}{rcl}
U_{\lambda,\mu}&\longrightarrow&
\CC[S_{\infty}\backslash B_-]_q[[\lambda^{-1},\mu^{-1}]]\\
u_{\lambda}&\longmapsto&u(\lambda)\\
u_{\mu}&\longmapsto&u(\mu).
\end{array}
\end{equation}
It can be shown that this morphism is injective.
There are unique coefficients $c_{\a,\b}^{k,l}$ 
satisfying 
\begin{equation}
   \forall\,k,l\in\NN,\quad (u_{\mu})^{k}(u_{\lambda})^{l}=
    \sum_{\a+\b=k+l}c_{\a,\b}^{k,l}\, (u_{\lambda})^{\a}(u_{\mu})^{\b}.
\end{equation}
We give below formulas which deal with cases $k$ or $l$ equal to $1$
or $2$. These formulas will be useful because the relations between 
$u_{\lambda}$ and $u_{\mu}$ are quadratic as well as the 
right side of (\ref{jmul}).
 \begin{prop}\label{propcab}
There are coefficients $c_{\a,\b},\, d_{\a,\b},\,
c_{\a,\b}^{(2)},\, d_{\a,\b}^{(2)}$ such that~:
\begin{align}
   \label{cab}
\forall\, n\in\NN,\quad u_{\mu}(u_{\lambda})^{n}&=\sum_{\a+\b=n+1}
c_{\a,\b}\, (u_{\lambda})^{\a}(u_{\mu})^{\b}\\
   \label{dab}
(u_{\mu})^{n}u_{\lambda}&=\sum_{\a+\b=n+1}
d_{\a,\b}\, (u_{\lambda})^{\a}(u_{\mu})^{\b}\\
   \label{cab2}
(u_{\mu})^{2}(u_{\lambda})^{n}&=\sum_{\a+\b=n+1}
c_{\a,\b}^{(2)}\, (u_{\lambda})^{\a}(u_{\mu})^{\b}\\
   \label{dab2}
(u_{\mu})^{n}(u_{\lambda})^{2}&=\sum_{\a+\b=n+1}
d_{\a,\b}^{(2)}\, (u_{\lambda})^{\a}(u_{\mu})^{\b}.
\end{align}
For all $\a\geq 0$ and $b\geq 1$, we have~:
\begin{align}
\label{fca0}
c_{\a,0}&=
\dfrac{(q^{\a-1}-1)\lambda^{-1}}{q^{\a-1}\lambda^{-1}-\mu^{-1}}\\
\label{fcab}
\forall\,\b\not=0,\quad c_{\a,\b}&=q^{\a-1}\Bigl[
\prod_{j=\a+1}^{\a+\b-1}(q^{j}-1)\Bigr]
\dfrac{(\lambda^{-1}-\mu^{-1})(\lambda^{-1}-q\mu^{-1})\mu^{-(\b-1)}}{
\prod_{j=\a-1}^{\a+\b-1}(q^{j}\lambda^{-1}-\mu^{-1})}.
\end{align}
Also, for all $\a\geq 0$ and $\b\geq 2$,
\begin{align}
\label{fca02}
c_{\a,0}^{(2)}&=\dfrac{(q^{\a-2}-1)(q^{\a-1}-1)\lambda^{-2}}{(
q^{\a-2}\lambda^{-1}-\mu^{-1})(q^{\a-1}\lambda^{-1}-\mu^{-1})}\\
\label{fca12}
c_{\a,1}^{(2)}&=q^{\a-2}(q^{\a-1}-1)[2]\dfrac{(\lambda^{-1}-\mu^{-1})
(\lambda^{-1}-q\mu^{-1})\lambda^{-1}}{(q^{\a-2}\lambda^{-1}-\mu^{-1})
(q^{\a-1}\lambda^{-1}-\mu^{-1})(q^{\a}\lambda^{-1}-\mu^{-1})}\\
\label{fcab2}
c_{\a,\b}^{(2)}&=q^{\a-2}\Bigl[
\prod_{j=\a+1}^{\a+\b-2}(q^{j}-1)\Bigr]\dfrac{(\lambda^{-1}-\mu^{-1})
(\lambda^{-1}-q\mu^{-1})\mu^{-(\b-2)}P_{\a,\b}^{(2)}(\lambda^{-1},\mu^{-1})}
{\prod_{j=\a-2}^{\a+\b-1}(q^{j}\lambda^{-1}-\mu^{-1})}
\end{align}
with~: 
\begin{align}
     \label{pab2}
       P_{\a,\b}^{(2)}(\lambda^{-1},\mu^{-1})&=q[\a+\b-1]
       (q\lambda^{-1}-\mu^{-1})(q^{\a-2}\lambda^{-1}-\mu^{-1})\\
     \notag
       &\phantom{=}-[\a]
       (\lambda^{-1}-q\mu^{-1})(q^{\a+\b-1}\lambda^{-1}-\mu^{-1}).
\end{align}
Coefficients $d_{\a,\b}$ (resp. $d_{\a,\b}^{(2)}$)
are obtained from $c_{\a,\b}$ (resp. $c_{\a,\b}^{(2)}$)
by~:
\begin{align}
\forall\,\a,\b,\quad
\lambda^{-(\b-1)}c_{\a,\b}&=\mu^{-(\b-1)}d_{\b,\a}\\
\lambda^{-(\b-2)}c_{\a,\b}^{(2)}&=\mu^{-(\b-2)}d_{\b,\a}^{(2)}.
\end{align}
\end{prop}

\begin{dem}
Coefficients $c_{\a,\b}^{(2)}$
and $d_{\a,\b}^{(2)}$ 
can be obtained by computation from $c_{\a,\b}$ 
and $d_{\a,\b}$. To prove  (\ref{fca0}) and (\ref{fcab}), we define
$c_{\a,0}$ and $c_{\a,\b}$ by these formulas, and we try to prove 
(\ref{cab}). For that, we express $u_{\lambda}$
and $u_{\mu}$ in terms of variables $v:=u_{\lambda}-u_{\mu}$
and $v':=\lambda^{-1}u_{\lambda}-\mu^{-1}u_{\mu}$.
These variables being $q$-commuting, we expand
the two expressions 
$\sum_{\a+\b=n+1}c_{\a,\b}\, (u_{\lambda})^{\a}(u_{\mu})^{\b}$
and $u_{\mu}(u_{\lambda})^{n}$ as a sum of terms in 
$v^{i}{v'}^{j}$. Then, we fix $i$ and $j$ and we want to
identify coefficients in $v^{i}{v'}^{j}$. This leads to
prove an equality between polynomials in 
$\dfrac{\lambda^{-1}}{\mu^{-1}}$ which reduces as a relation between
$q$-integers. In the same manner, we prove (\ref{dab}).
\end{dem}
\subsection{Proof of Proposition \ref{proex}}
First, we define derivations $H_n$ on the free algebra ${\cal A}$
generated by $u_i,\, m_i,\, i>0,\,\SPB$ and $\SMB$.
To prove that the $H_n$ give derivations on 
$\CC[H_-\backslash B_-N_+]_q$, we have to prove several relations.
The most complicated one is
\begin{equation}\label{relhnulum}
\pi\circ H(\nu).\bigl({\text{Relation between }}u(\lambda){\text{ and }}
u(\mu)\bigr)=0
\end{equation}
where $\pi$ denotes the projection of ${\cal A}$
onto $\CC[H_-\backslash B_-N_+]_q$.
To prove (\ref{relhnulum}), we decompose 
$H(\nu).u(\lambda)$ and $H(\nu).u(\mu)$
according to the relation~:
\begin{multline}
\label{devhmul}
H(\mu).u(\lambda)=\mu^{-1}(u(\mu)-u(\lambda))
+\sum_{k=0}^{+\infty}(-1)^{k+1}q^{\{ k+1\}}\times\\
\times\bigl(
\lambda^{-1}u(\lambda)u(\mu)^{k}u(\lambda)+\mu^{-1}
u(\mu)^{k+2}-\mu^{-1}u(\mu)^{k+1}u(\lambda)-
\mu^{-1}u(\lambda)u(\mu)^{k+1}\bigr)
m(\mu)^{k+1}\mu^{-(k+1)}
\end{multline}
which comes from (\ref{jmul}) and (\ref{perum}) by decomposing $v(\mu)$
into a generating series in $u(\mu)^k m(\mu)^k \mu^{-k}$.

\noindent
Thus, the left side of (\ref{relhnulum})
is of the form 
$\sum_{k=0}^{+\infty}
P_k(u(\lambda),u(\mu),u(\nu))m(\nu)^{k}\nu^{-k}$
where $P_k(u(\lambda),u(\mu),u(\nu))$ is a polynomial
in non-commutative variables $u(\lambda),\,u(\mu)$ and $u(\nu)$.

\noindent
Let $k\in\NN$. According to the fact that the relations between
$u(\lambda)$ and $u(\mu)$ are quadratric and that the only
terms in $u$ which appear in (\ref{jmul}) are also quadratic, 
it is possible to reorganize the terms of polynomial 
$P_{k}(u(\lambda),u(\mu),u(\nu))$ by using
Proposition \ref{propcab} and the morphism defined in
(\ref{morspeu}) so as to obtain a sum of monomials of the form
$(u(\lambda))^{\a}(u(\mu))^{\b}(u(\nu))^{\gamma}$.
Then, we fix $\a,\b,\gamma$ and we show that the coefficient
of $(u(\lambda))^{\a}(u(\mu))^{\b}(u(\nu))^{\gamma}$ in $P_k$
is equal to $0$. Thus, $P_k=0$, and
(\ref{relhnulum}) is true.
In the same way, we prove all other relations.
Thus, $H_n$ exists. To prove the commutativity,
we deduce from formulas
(\ref{jmul}) and (\ref{jmum}) that
\begin{align}
\label{imvl}
H(\mu)(v(\lambda))&=\dfrac{1}{\lambda^{-1}-\mu^{-1}}
\bigl(\mu^{-1}v(\lambda)-\lambda^{-1}v(\mu)\bigr)
+\dfrac{1}{\lambda^{-1}-\mu^{-1}}
v(\lambda)\bigl(\mu^{-1}u(\mu)-\lambda^{-1}u(\lambda)\bigr)v(\mu)\\
\notag
&+\dfrac{1}{\lambda^{-1}-\mu^{-1}}
v(\mu)\bigl(\mu^{-1}u(\mu)-\lambda^{-1}u(\lambda)\bigr)v(\lambda).
\end{align}
Then, the computation shows that 
$H(\mu)\circ H(\nu)=H(\nu)\circ H(\mu)$
(it is not necessary for that to decompose in generating series).
On the other hand, with the help of the formulas coming from
Corrolary \ref{stru},
the computation shows that
$H_{\mu}\circ f_{\pm}=f_{\pm}\circ H_{\mu}$.
This completes the proof of Proposition \ref{proex}, 
and also of Proposition \ref{exih}.
\subsection{End of the proof of Theorem \ref{thc}}
We are going to prove that for any integer $n$,
$\ad(I_n)\circ\DSB_q=\DSB_q\circ H_n$.
In other words, $\DSB_q$ set up a Drinfeld-Sokolov
correspondence for
the extended phase space ${\bar A}_q$
and the quantum homogeneous space 
$\CC[H_-\backslash B_-]_q$. Thanks to Lemma \ref{plon},
the result will follow.
To simplify, we note $U_q=\CC[H_-\backslash B_-]_q$ and
${\bar U}_q=\CC[H_-\backslash B_-N_+]_q$.
In the following of the article, we shall identify the elements
of ${\bar U}_q$ with their images in ${\bar A}_q$
by the algebra monomorphism $\DSB_q$.
The following lemma will be useful.

\begin{lem}\label{uniprod}
Let $D_1$ and $D_2$ be two derivations
defined on $A_q$ (resp. ${\bar A}_q$) such that 
$D_1(x)=D_2(x)$ for all $x\in U_q$ (resp. ${\bar U}_q$).
Then, $D_1=D_2$.
\end{lem}
 
\begin{dem}
Let us denote by ${\cal C}$ the subalgebra of $A_q$
(resp. ${\bar A}_q$) of all elements $x$ such that
$D_1(x)=D_2(x)$. If $x\in{\cal C}$, is invertible in $A_q$
(resp. ${\bar A}_q$) then $x^{-1}\in{\cal C}$.
By induction, using the explicit forms of $u_n$ and $m_n$ given in 
(\ref{foru}) et (\ref{defmiq}), we show that
for all $n\in\ZZ$, $x_n^{\pm 1},y_n^{\pm 1}\in{\cal C}$.
For example, $y_0^{-1}\in{\cal C}$,
$x_0^{-1}=q^{-1}u_2 y_0^2\in{\cal C}$.
Hence, we get the result.
\end{dem}

Let $n\in\NN^*$. Then $H_n$ is a derivation on 
${\bar U}_q\subset{\bar A}_q$. We have to prove that 
$H_n$ extends as a derivation on ${\bar A}_q$ and that
$H_n=\ad(I_n)$.

\subsubsection{First, let us assume that $H_n$ 
has an extension to $A_q$}
Then, $H_n$ has also an extension on ${\bar A}_q$.
Moreover, thanks to the definition of $H_n$ 
together with the relations (\ref{fracu}) and (\ref{fracm})
which give expressions for $u(\lambda)$ 
and $m(\lambda)$ as quantum continued fractions,
we show that the image of ${\bar U}_q$ by 
$T^{-\frac{1}{2}}$ is generated by ${\bar U}_q$
and $u_1^{-1}$ and that 
$H_n\circ T^{-\frac{1}{2}}(x)=T^{-\frac{1}{2}}\circ H_n(x)$
for all $x\in{\bar U}_q$. Lemma \ref{uniprod}
ensures that this relation is also true on ${\bar A}_q$.
On the other hand, computation shows that
for all $x\in U_q,\, H_n\circ\varphi(x)=\varphi\circ H_n(x)$ 
where $\varphi$ is the involution on $A_q$ defined in (\ref{definfi}).
Using again Lemma \ref{uniprod}, this last equality
extends on $A_q$. So, by using the relation
$\varphi\circ T^{ \frac{1}{2}}=T^{-\frac{1}{2}}\circ\varphi$,
we deduce that $H_n\circ T^{ \frac{1}{2}}(x)= T^{ \frac{1}{2}} \circ
H_n(x)$ for all $x\in A_q$. It can be shown 
that this relation is also true for $x=\SPMB$.
It follows that $H_n\circ T^{\pm\frac{1}{2}}=
T^{\pm\frac{1}{2}}\circ H_n$ on ${\bar A}_q$. 
On the other hand, by virtue of Proposition \ref{proex}, we have 
$H_n\circ f_{\pm} (y_0^{-1})=f_{\pm}\circ H_n (y_0^{-1})$
for $u_1=y_0^{-1}$. {}From the equality 
$f_{\pm}\circ T^{-\frac{1}{2}}=T^{-\frac{1}{2}}\circ f_{\mp}$,
we deduce that $H_n\circ f_{\pm}(x)=f_{\pm}\circ H_n(x)$
for all $x\in B_q$, where $B_q$ denotes the subalgebra
of $A_q$ generated by $x_i^{-1}$ and $y_i^{-1}$, for $i\in\ZZ$.
Then, the $U_q\bm$-module-algebra structure of ${\bar A}_q$
implies that $H_n\circ f_{\pm}(x)=f_{\pm}\circ H_n(x)$
for all $x\in A_q$ and also on ${\bar A}_q$
for this equality is also true if
$x=\SEPB$, $\varepsilon\in\{ +,-\}$.
Thus, $H_n\in\Der_{U_q\bm}({\bar A}_q)$ and thanks to Proposition
\ref{onv}, we see that $H_n$ is a linear combination of $\ad(I_k)$.
But, the same proposition also shows that there is a gradation
$\deg$ on $\Der_{U_q\bm}({\bar A}_q)$ given by $\deg \delta=n$
for $\delta\in\Der_{U_q\bm}({\bar A}_q)$ if
there is a homogeneous element $\a\in {\bar A}_q$ with respect
to the principal gradation $\deg_p$ on ${\bar A}_q$
defined on $A_q$ in subsection \ref{defgra}
and extended on ${\bar A}_q$ by $\deg_p\SPM=1$,
such that $\delta(\a)$ is also homogeneous with respect to $\deg_p$
and $\deg_p\delta(\a)=n+\deg_p(\a)$.
Note that if $\delta$ is homogeneous,
this last property occurs not only for one special homogeneous
element $\a$ with respect to the gradation $\deg_p$
but also for all homogeneous elements $x$ in ${\bar A}_q$
with respect to $\deg_p$.
Now, from the definition of $H_n$, it is easy to see that
$\deg H_n=\deg I_n$. Then, the computation of 
both $H_n(y_0^{-1})$ and $\ad(I_n)(y_0^{-1})$
on the basis element 
$y_{-k}^{-1}x_{-k+1}^{-2}y_{-k+1}^{-2}\ldots
x_0^{-2}y_0^{-2}$ or $x_{-k}^{-1}y_{-k}^{-2}\ldots
x_0^{-2}y_0^{-2}$ 
of the basis $\prod x_i^{\alpha_i}y_i^{\beta_i}$
(with $k$ such that $n=2k$ or $n=2k+1$ according to the parity of $n$)
shows that $H_n=\ad(I_n)$.
Thus, to conclude, it suffices to extend $H_n$ on ${\bar A}_q$.
\subsubsection{Proof of the existence of an extension}
In the classical case,
we use the fact that the classical limit $U_{\cl}$ of $U_q$
possesses the same field of fraction $K$ as $A_{\cl}$
to extend $H_{n,\cl}$ in a derivation of $K$.
The extension is unique. So, the relation 
$H_{n,\cl}\circ T^{-\frac{1}{2}}=T^{-\frac{1}{2}}\circ H_{n,\cl}$
true on $A_{\cl}$ is also true on $K$.
The same argument as above with the involution $\varphi$
shows that 
$H_{n,\cl}\circ T^{\pm\frac{1}{2}}=T^{\pm\frac{1}{2}}\circ H_{n,\cl}$.
But, $H_{n,\cl}(y_0^{-1})=H_{n,\cl}(u_1)\in U_{\cl}\subset A_{\cl}$.
Thus, $A_{\cl}$ is invariant by $H_n$. Hence, we get the result.

In the classical case, it is a little bit more complicated.
The problem comes from the fact that it is not obvious that 
non-zero elements of $A_q$ as well as of $\CC[H_-\backslash B_-]_q$
satisfy Ore conditions. Nevertheless, we show that this is true when $q$
is a {\it formal} variable ``close to $1$''. For that,
we develop the notion of extended Ore conditions.

\begin{Def}
{\it{Let $A$ be an algebra without any zero divisors
over a field $k$ and $(A[[t]],*)$ a formal deformation
of the multiplication on $A$. For all $n$, we note by $\pi_n$ the natural 
projection of $A[[t]]$ on $A_n:=A[[t]]/(t^n)$. A multiplicative set $S$
in $A[[t]]$ is said to satisfy the extended Ore conditions
if $\pi_n(S)$ satisfy the Ore conditions in $A_n$ equipped with
the natural non-commutative product induced by $*$.}}
\end{Def}

If $S$ is a multiplicative set which satisfies the extended Ore conditions
in $A[[t]]$, then there are natural morphisms~:
$(A_n)_{S_n}\longrightarrow (A_p)_{S_p}$ pour $n>p$. We note
$A[[t]]_S$ the projective limit of $(A_n)_{S_n}$.

\begin{ex}

1. Let $B=\CC[X_i^{-1},Y_i^{-1},\, i\in\ZZ]$, 
and $K$ be the field of fractions of $B$.
Let us consider the isomorphism of free
modules~:
$$
\begin{array}{rcl}
B[[q-1]]&\longrightarrow&B'_q\\
\prod_{i=1}^{\infty}X_i^{-\alpha_i}\prod_{j=1}^{\infty}Y_j^{-\beta_j}&
\longmapsto&
\prod_{i=1}^{\infty}x_i^{-\alpha_i}\prod_{j=1}^{\infty}y_j^{-\beta_j}
\end{array}
$$
where $(\alpha_i)$ and $(\beta_j)$
are two almost zero sequences in
$\NN^{\ZZ}$ and $B'_q$ is the $(q-1)$-adic completion of the
subalgebra
$B_q$ aforementioned. 
This isomorphism leads to a formal deformation
$*'$ of the multiplication on $B[[q-1]]$. 
We show that the multiplicative set 
$S'$ of all elements 
non divisible by $q-1$ satisfy the extended
Ore conditions.
Moreover, $*'$ is a star-product, i.e., 
there are $B_n$ bidifferential operators such that
$*'=\sum B_n (q-1)^n$.
Hence, we get a non-commutative structure on 
$(K[[q-1]],*')$ which contains
$(B[[q-1]],*')$ 
and it can be proved that
$(B[[q-1]]_{S'},*')\simeq (K[[q-1]],*')$.

2. One can also define a structure of non-commutative algebra 
$(\CC[U_i,M_i,i>0][[q-1]],*)$ thanks to the isomorphism of
free modules~:
$$
\begin{array}{rcl}
\CC[U_i,M_i,i>0][[q-1]]&\longrightarrow&U'_q\\
\prod_{i=1}^{\infty}U_i^{\alpha_i}M_i^{\beta_i}&\longmapsto&
\prod_{i=1}^{\infty}u_i^{\alpha_i}m_i^{\beta_i}
\end{array}
$$
where $U'_q$ denotes the $(q-1)$-adic completion of
$U_q$. Contrary to the previous case, it is not so easy to
check that the non-commutative product $*$ is a star product.
So, there is {\it a priori} no reason for
$(\CC(U_i,M_i,i>0)[[q-1]],*)$ to exist.
\end{ex}
However, we prove the following result. 

\begin{lem}\label{orelem}
For all $x\in A_q$, there is $\omega$ a monomial in
$x_i^{-1},\, y_i^{-1}$ such that $\omega\in U_q$ and 
$\omega x\in U_q$.
\end{lem}

\begin{dem}
It suffices to prove the lemma for $x=x_i^{-1}$ or $x=y_i^{-1}$
with $i\in\ZZ$. Thanks to the symmetry relation between $u_n$
and $m_n$ i.e., $m_n=\varphi(u_n)$ and 
$u_n\in\CC[x_j^{-1},y_j^{-1},\, j<0]_q$,
we can assume that $i<0$.
Then, we prove the result by induction.
For example, $1.y_0^{-1}=u_1$, $q^{-1}y_0^{-2}x_0^{-1}=u_2$
and $y_0^{-1}=u_1$, etc. 
\end{dem}

Lemma \ref{orelem} implies that the multiplicative set $S$
of all elements in $\CC[U_i,M_i,i>0][[q-1]]$
which are non-divisible by $q-1$ satisfy the extended
Ore conditions and that we have an isomorphism
$(\CC[U_i,M_i,i>0][[q-1]]_S,*)\simeq (B[[q-1]]_{S'},*')$.
Now, to conclude, we say that the $(q-1)$-adic completion of $H_n$
defined on $\CC[U_i,M_i,i>0][[q-1]]$ 
has naturally an extension on $(K[[q-1]],*')$.
Moreover, this extension is unique. The same arguments as above
with the involution $\varphi$ and the half-translation automorphism
$T^{\pm\frac{1}{2}}$ show that $H_n$
and $T^{\pm\frac{1}{2}}$ commute on $K[[q-1]]$.
But $H_n(y_0^{-1})\subset A_q$. It follows that
$A_q$ is invariant by $H_n$.

\section{Conclusion and outlooks}
First of all, it would be interesting to see whether if it is possible
to extend our result to the more general case of an arbitrary
non-twisted Lie algebra and to study other possible
models of discretization proposed by Enriquez and Feigin
\cite{ENR}.
It would be also interesting to study in details
the case when $q$ is a root of the unity.

\subsection{Affine Poisson homogeneous space}
Theorem \ref{thc} leaves us with the feeling that there is a 
general Drinfeld-Sokolov correspondence for the discrete Toda theory
and that the homogeneous space of the correspondence
is a Poisson homogeneous space equipped with 
a Poisson structure induced by a Poisson bivector $\pi$
of the form $\pi=r^L-{r'}^R$, where $r$ and $r'$ are two $r$-matrices
such that their Schouten bracket $[r,r]$ and $[r',r']$
are equal and invariant by the adjoint action of the Lie group
$G$ on $\bigwedge^3\,\text{Lie}(G)$ and where $r^L$
(resp. $r^R$) is the left (resp. right)
translation of $r$ (resp. $r'$) on $G$. A group $G$
endowed with such a Poisson structure is a particular case 
of an affine Poisson homogeneous space (APHS),
according to the terminology introduced by Dazord and Sondaz
\cite{DaSo} (see also \cite{LU} and \cite{KOS}).
By definition, an APHS is a Poisson manifold which is a
principal homogeneous space under the action of a Poisson-Lie group
and that if it is the case, then there are on $G$ two commuting actions
by Poisson-Lie groups.
\subsection{Links with Parmentier's work}
Our method of quantizing the Poisson manifold 
$H_-\backslash B_-$ equipped with the Poisson structure induced by the field of bivectors 
$P_{\infty}=r^L-r_{\infty}^R$ 
(which is truly a quotient of an APHS) lays on the study of the classical case
and on the fact that the phase space of the discrete sine-Gordon system 
had a natural quantization. But, it is perhaps
not the easiest way to quantize $(H_-\backslash B_-,P_{\infty})$. Indeed, 
in the case we dealt with, $r$ denotes the standard $r$-matrix and $r_{\infty}$ denotes the
$r$-matrix corresponding to Drinfeld new realizations.
So, according to Parmentier works, to get a quantization of the APHS $G$
with Poisson structure given by the field of bivectors $P_{\infty}$,
it suffices to have a twist relying the two Hopf algebra structures
$(U_q{\mathfrak{g}},\Delta)$ and
$(U_q{\mathfrak{g}},\Delta_{nr})$ where, in the first case,
the comultiplication is the ``Drinfeld-Jimbo'' comultiplication, and in the second case,
it is the one corresponding to the Drinfeld new realizations.
But such a twist appears in the paper \cite{KT}. 
We need to apply this method and to investigate further 
about how derivations $H_n$ appear in it.
We plan to study this question elsewhere.
\subsection{The continuous case}
It would be also interesting to see whether it is possible to deduce from our results
solutions to problems of continuous Toda theory, 
to compute explicitly integrals of motion,
to quantize in terms of Vertex Operator Algebra the 
Vertex Poisson Algebra shown by Enriquez and Frenkel on homogeneous spaces
\cite{EFRE}, and to obtain a quantum version of Drinfeld-Sokolov correspondence
in terms of V.O.A. in the continuous case.
\subsection{The Drinfled-Sokolov reduction}
Finally, we indicate that there exists another correspondence, close to the one
we discussed here, which is called the Drinfeld-Sokolov {\it reduction}
\cite{DS}. This correspondence allows to construct ${\cal W}$-algebras
from Kac-Moody algebras. There is a Poisson isomorphism between the manifold
of scalar differential operators of order $n$ with the second
Gelfand-Dickey bracket on one hand, and the manifold of matrix differential
operators of order $1$ viewed as a subspace of 
${\widehat{\mathfrak{sl}_n}}^*$,
with Kirillov-Kostant bracket on the other hand.
The quantization of this correspondence is studied in
\cite{FFA}. A $q$-deformed version of this correspondence, in which manifolds of differential
operators are replaced by manifolds of $q$-difference operators is proposed in
\cite{FRS} and \cite{SS}.
Quantization of this correspondence leads to $q$-deformed ${\cal W}$-algebra.

\bigskip

{\bf Acknowledgement.}
I am very grateful to my teacher Professor Benjamin
Enriquez for his help and support during the preparation 
of this paper. I would also like to thank Professor
Anthony Joseph and Professor Joseph Bernstein for the 
hospitality of the Weizmann institute and the Tel-Aviv 
university.

{}

\end{document}